# The Beta Generalized Marshall-Olkin-G Family of Distributions

Laba Handique and Subrata Chakraborty*
Department of Statistics, Dibrugarh University
Dibrugarh-786004, India
*Corresponding Author. Email: subrata_stats@dibru.ac.in

**Abstract**

In this paper we propose a new family of distribution considering Generalized Marshal-Olkin distribution as the base line distribution in the Beta-*G* family of Construction. The new family includes Beta-*G* (Eugene et al. 2002 and Jones, 2004) and *GMOE* (Jayakumar and Mathew, 2008) families as particular cases. Probability density function (pdf) and the cumulative distribution function (cdf) are expressed as mixture of the Marshal-Olkin (Marshal and Olkin, 1997) distribution. Series expansions of pdf of the order statistics are also obtained. Moments, moment generating function, Rényi entropies, quantile power series, random sample generation and asymptotes are also investigated. Parameter estimation by method of maximum likelihood and method of moment are also presented. Finally proposed model is compared to the Generalized Marshall-Olkin Kumaraswamy extended family (Handique and Chakraborty, 2015) by considering three data fitting examples with real life data sets.

## 1. Introduction

Here we briefly introduce the Beta-*G* (Eugene et al. 2002 and Jones, 2004) and Generalized Marshall-Olkin family (Jayakumar and Mathew, 2008) of distributions**.**

### 1.1 Some formulas and notations

Here first we list some formulas to be used in the subsequent sections of this article.

If $T$ is a continuous random variable with pdf, $f(t)$ and cdf $F(t) = P[T \le t]$, then its

Survival function (sf): $\bar{F}(t) = P[T > t] = 1 - F(t)$,

Hazard rate function (hrf): $h(t) = f(t)/\bar{F}(t)$,

Reverse hazard rate function (rhrf): $r(t) = f(t)/F(t)$,

Cumulative hazard rate function (chrf): $H(t) = -\log\,[\overline{F}(t)]$,

$(p,q,r)^{\text{th}}$ Power Weighted Moment (PWM): $\Gamma_{p,q,r} = \int_{-\infty}^{\infty} t^p\,[F(t)]^q\,[1-F(t)]^r\, f(t)\,dt$,

Rényi entropy: $I_R(\delta) = (1-\delta)^{-1}\log\left(\int_{-\infty}^{\infty} f(t)^{\delta}\,dt\right)$.

**1.2 Beta-*G* family of distributions**

The cdf of beta-*G* (Eugene et al. 2002 and Jones 2004) family of distribution is

$$F^{BG}(t) = \frac{1}{B(m,n)}\int_{o}^{F(t)} v^{m-1}(1-v)^{n-1}\,dv = \frac{B_{F(t)}(m,n)}{B(m,n)} = I_{F(t)}(m,n) \qquad (1)$$

Where $I_t(m,n) = B(m,n)^{-1}\int_0^t x^{m-1}(1-x)^{n-1}dx$ denotes the incomplete beta function ratio.

The pdf corresponding to (1) is

$$f^{BG}(t) = \frac{1}{B(m,n)} f(t)F(t)^{m-1}[1-F(t)]^{n-1} = \frac{1}{B(m,n)} f(t)F(t)^{m-1}\overline{F}(t)^{n-1} \qquad (2)$$

Where $f(t) = dF(t)/dt$.

sf: $\overline{F}^{BG}(t) = 1 - I_{F(t)}(m,n) = \dfrac{B(m,n) - B_{F(t)}(m,n)}{B(m,n)}$

hrf: $h^{BG}(t) = f^{BG}(t)/\overline{F}^{BG}(t) = \dfrac{f(t)F(t)^{m-1}\overline{F}(t)^{n-1}}{B(m,n) - B_{F(t)}(m,n)}$

rhrf: $r^{BG}(t) = f^{BG}(t)/F^{BG}(t) = \dfrac{f(t)F(t)^{m-1}\overline{F}(t)^{n-1}}{B_{F(t)}(m,n)}$

chrf: $H(t) = -\log\left[\dfrac{B(m,n) - B_{F(t)}(m,n)}{B(m,n)}\right]$

Some of the well known beta generated families are the Beta-generated (beta-*G*) family (Eugene *et al.*, 2002; Jones 2004), beta extended *G* family (Cordeiro *et al.*, 2012), Kumaraswamy beta generalized family (Pescim *et al.*, 2012), beta generalized weibull distribution (Singla et al., 2012), beta generalized Rayleigh distribution (Cordeiro *et al.*, 2013), beta extended half normal distribution (Cordeiro *et al.*, 2014), beta log-logistic distribution (Lemonte, 2014) beta generalized inverse Weibull distribution (Baharith et al., 2014), beta Marshall-Olkin family of distribution (Alizadeh *et al.*, 2015) and beta generated Kumaraswamy-G family of distribution (Handique and Chakraborty 2016a) among others.

**1.3 Generalized Marshall-Olkin Extended (GMOE) family of distributions**

Jayakumar and Mathew (2008) proposed a generalization of the Marshall and Olkin (1997) family of distributions by using the Lehman second alternative (Lehmann 1953) to obtain the sf $\overline{F}^{GMO}(t)$ of the GMOE family of distributions by exponentiation the sf of MOE family of distributions as

$$\overline{F}^{GMO}(t) = \left[\frac{\alpha\,\overline{G}(t)}{1-\overline{\alpha}\,\overline{G}(t)}\right]^{\theta}, -\infty < t < \infty; 0 < \alpha < \infty; 0 < \theta < \infty \quad (3)$$

where $-\infty < t < \infty$, $\alpha > 0\,(\overline{\alpha} = 1-\alpha)$ and $\theta > 0$ is an additional shape parameter. When $\theta = 1$, $\overline{F}^{GMO}(t) = \overline{F}^{MO}(t)$ and for $\alpha = \theta = 1$, $\overline{F}^{GMO}(t) = \overline{F}(t)$. The cdf and pdf of the GMOE distribution are respectively

$$F^{GMO}(t) = 1 - \left[\frac{\alpha\,\overline{G}(t)}{1-\overline{\alpha}\,\overline{G}(t)}\right]^{\theta} \quad (4)$$

and $$f^{GMO}(t) = \theta\left[\frac{\alpha\,\overline{G}(t)}{1-\overline{\alpha}\,\overline{G}(t)}\right]^{\theta-1}\left\{\frac{\alpha\,g(t)}{[1-\overline{\alpha}\,\overline{G}(t)]^{2}}\right\} = \frac{\theta\alpha^{\theta}\,g(t)\overline{G}(t)^{\theta-1}}{[1-\overline{\alpha}\,\overline{G}(t)]^{\theta+1}} \quad (5)$$

Reliability measures like the hrf, rhrf and chrf associated with (1) are

$$h^{GMO}(t) = \frac{f^{GMO}(t)}{\overline{F}^{GMO}(t)} = \frac{\theta\alpha^{\theta}\,g(t)\overline{G}(t)^{\theta-1}}{[1-\overline{\alpha}\,\overline{G}(t)]^{\theta+1}} \Bigg/ \left[\frac{\alpha\,\overline{G}(t)}{1-\overline{\alpha}\,\overline{G}(t)}\right]^{\theta}$$

$$= \theta\frac{g(t)}{\overline{G}(t)}\frac{1}{1-\overline{\alpha}\,\overline{G}(t)} = \theta\,\frac{h(t)}{1-\overline{\alpha}\,\overline{G}(t)}$$

$$r^{GMO}(t) = \frac{f^{GMO}(t)}{F^{GMO}(t)} = \frac{\theta\alpha^{\theta}\,g(t)\overline{G}(t)^{\theta-1}}{[1-\overline{\alpha}\,\overline{G}(t)]^{\theta+1}} \Bigg/ 1 - \left[\frac{\alpha\,\overline{G}(t)}{1-\overline{\alpha}\,\overline{G}(t)}\right]^{\theta}$$

$$= \frac{\theta\alpha^{\theta}\,g(t)\overline{G}(t)^{\theta-1}}{[1-\overline{\alpha}\,\overline{G}(t)]\;[\{1-\overline{\alpha}\,\overline{G}(t)\}^{\theta} - \alpha^{\theta}\,\overline{G}(t)^{\theta}\,]} = \frac{\theta\alpha^{\theta}\,g(t)\overline{G}(t)^{\theta-1}}{[1-\overline{\alpha}\,\overline{G}(t)]^{\theta+1} - \alpha^{\theta}\,\overline{G}(t)^{\theta}[1-\overline{\alpha}\,\overline{G}(t)]}$$

$$H^{GMO}(t) = -\log\left[\frac{\alpha\,\overline{G}(t)}{1-\overline{\alpha}\,\overline{G}(t)}\right]^{\theta} = -\theta\log\left[\frac{\alpha\,\overline{G}(t)}{1-\overline{\alpha}\,\overline{G}(t)}\right]$$

Where $g(t)$, $G(t)$, $\overline{G}(t)$ and $h(t)$ are respectively the pdf, cdf, sf and hrf of the baseline distribution. We denote the family of distribution with pdf (3) as GMOE$(\alpha,\theta,a,b)$ which for $\theta = 1$, reduces to MOE$(\alpha, a,b)$ .

Some of the notable distributions derived Marshall-Olkin Extended exponential distribution (Marshall and Olkin, 1997), Marshall-Olkin Extended uniform distribution (Krishna, 2011; Jose and Krishna 2011), Marshall-Olkin Extended power log normal distribution (Gui, 2013a), Marshall-Olkin Extended log logistic distribution (Gui, 2013b), Marshall-Olkin Extended Esscher transformed Laplace distribution (George and George, 2013), Marshall-Olkin Kumaraswamy-*G* familuy of distribution (Handique and Chakraborty, 2015a), Generalized Marshall-Olkin Kumaraswamy-*G* family of distribution (Handique and Chakraborty 2015b) and Kumaraswamy Generalized Marshall-Olkin family of distribution (Handique and Chakraborty 2016b).

In this article we propose a family of Beta generated distribution by considering the Generalized Marshall-Olkin family (Jayakumar and Mathew, 2008) as the base line distribution in the Beta-*G* family (Eugene et al. 2002 and Jones, 2004). This new family referred to as the Beta Generalized Marshall-Olkin family of distribution is investigated for some its general properties. The rest of this article is organized in seven sections. In section 2 the new family is introduced along with its physical basis. Important special cases of the family along with their shape and main reliability characteristics are presented in the next section. In section 4 we discuss some general results of the proposed family. Different methods of estimation of parameters along with three comparative data modelling applications are presented in section 5. The article ends with a conclusion in section 6 followed by an appendix to derive asymptotic confidence bounds.

## 2. New Generalization: Beta Generalized Marshall-Olkin-G ($BGMO-G$) family of distributions

Here we propose a new Beta extended family by considering the cdf and pdf of $GMO$ (Jayakumar and Mathew, 2008) distribution in (4) and (5) as the $f(t)$ and $F(t)$ respectively in the Beta formulation in (2) and call it $BGMO-G$ distribution. The pdf of $BGMO-G$ is given by

$$f^{BGMOG}(t)=\frac{1}{B(m,n)}\frac{\theta\alpha^{\theta}g(t)\bar{G}(t)^{\theta-1}}{[1-\bar{\alpha}\bar{G}(t)]^{\theta+1}}\left[1-\left[\frac{\alpha\bar{G}(t)}{1-\bar{\alpha}\bar{G}(t)}\right]^{\theta}\right]^{m-1}\left[\left[\frac{\alpha\bar{G}(t)}{1-\bar{\alpha}\bar{G}(t)}\right]^{\theta}\right]^{n-1} \quad (6)$$

$$,\ 0<t<\infty,\ 0<a,b<\infty, m>0, n>0$$

Similarly substituting from equation (4) in (1) we get the cdf of $BGMO-G$ respectively as

$$F^{BGMOG}(t)=I_{1-\left[\frac{\alpha\bar{G}(t)}{1-\bar{\alpha}\bar{G}(t)}\right]^{\theta}}(m,n) \quad (7)$$

The sf, hrf, rhrf and chrf of $BGMO-G$ distribution are respectively given by

sf: $$\bar{F}^{BGMOG}(t)=1-I_{1-\left[\frac{\alpha\bar{G}(t)}{1-\bar{\alpha}\bar{G}(t)}\right]^{\theta}}(m,n)$$

hrf :

$$h^{BGMOG}(t)=\frac{1}{B(m,n)}\frac{\theta\alpha^{\theta}g(t)\bar{G}(t)^{\theta-1}}{[1-\bar{\alpha}\bar{G}(t)]^{\theta+1}\,[1-I_{1-\left[\frac{\alpha\bar{G}(t)}{1-\bar{\alpha}\bar{G}(t)}\right]^{\theta}}(m,n)]}$$
$$\times\left[1-\left[\frac{\alpha\bar{G}(t)}{1-\bar{\alpha}\bar{G}(t)}\right]^{\theta}\right]^{m-1}\left[\left[\frac{\alpha\bar{G}(t)}{1-\bar{\alpha}\bar{G}(t)}\right]^{\theta}\right]^{n-1} \quad (8)$$

rhrf :

$$r^{BGMOG}(t)=\frac{1}{B(m,n)}\frac{\theta\alpha^{\theta}g(t)\overline{G}(t)^{\theta-1}}{[1-\overline{\alpha}\,\overline{G}(t)]^{\theta+1}\,[I_{1-\left[\frac{\alpha\overline{G}(t)}{1-\overline{\alpha}\,\overline{G}(t)}\right]^{\theta}}(m,n)]}$$

$$\times\left[1-\left[\frac{\alpha\overline{G}(t)}{1-\overline{\alpha}\,\overline{G}(t)}\right]^{\theta}\right]^{m-1}\left[\left[\frac{\alpha\overline{G}(t)}{1-\overline{\alpha}\,\overline{G}(t)}\right]^{\theta}\right]^{n-1} \quad (9)$$

chrf: $H^{BGMOG}(t)=-\log[\,1-I_{1-\left[\frac{\alpha\overline{G}(t)}{1-\overline{\alpha}\,\overline{G}(t)}\right]^{\theta}}(m,n)]$

Remark: The $\text{BGMO}-G(\text{m},\text{n},\alpha,\theta)$ reduces to

(i) $\text{BMO}(\text{m},\text{n},\alpha)$ (Alizadeh *et al.*, 2015) for $\theta=1$; (ii) $\text{GMO}(\alpha,\theta)$,(Jayakumar and Mathew, 2008) if $m=n=1$; (iii) $\text{MO}(\alpha)$ (Marshall and Olkin, 1997) when $m=n=\theta=1$; and (iv) $\text{B}(\text{m},\text{n})$ (Eugene *et al.*, 2002; Jones 2004) for $\alpha=\theta=1$.

**2.1 Genesis of the distribution**

If $m$ and $n$ are both integers, then the probability distribution takes the same form as the order statistics of the random variable $T$.

**Proof:** Let $T_1,T_2...,T_{m+n-1}$ be a sequence of *i.i.d.* random variables with cdf $1-[1-G(t)^{a}]^{b}$. Then the pdf of the $m^{th}$ order statistics $T_{(m)}$ is given by

$$=\frac{(m+n-1)!}{(m-1)!\,[(m+n-1)-m]!}[1-\{\alpha\overline{G}(t)/1-\overline{\alpha}\,\overline{G}(t)\}^{\theta}]^{m-1}[\{\alpha\overline{G}(t)/1-\overline{\alpha}\,\overline{G}(t)\}^{\theta}]^{(m+n-1)-m}$$

$$\times\theta\alpha^{\theta}g(t)\overline{G}(t)^{\theta-1}/[1-\overline{\alpha}\,\overline{G}(t)]^{\theta+1}$$

$$=\frac{\Gamma(m+n)}{\Gamma(m)\,\Gamma(n)}[1-\{\alpha\overline{G}(t)/1-\overline{\alpha}\,\overline{G}(t)\}^{\theta}]^{m-1}[\{\alpha\overline{G}(t)/1-\overline{\alpha}\,\overline{G}(t)\}^{\theta}]^{n-1}$$

$$\times\theta\alpha^{\theta}g(t)\overline{G}(t)^{\theta-1}/[1-\overline{\alpha}\,\overline{G}(t)]^{\theta+1}$$

$$=\frac{1}{B(m,n)}[1-\{\alpha\overline{G}(t)/1-\overline{\alpha}\,\overline{G}(t)\}^{\theta}]^{m-1}[\{\alpha\overline{G}(t)/1-\overline{\alpha}\,\overline{G}(t)\}^{\theta}]^{n-1}$$

$$\times\theta\alpha^{\theta}g(t)\overline{G}(t)^{\theta-1}/[1-\overline{\alpha}\,\overline{G}(t)]^{\theta+1}$$

**2.2 Shape of the density function**

Here we have plotted the pdf of the $BGMO-G$ for some choices of the distribution $G$ and parameter values to study the variety of shapes assumed by the family.

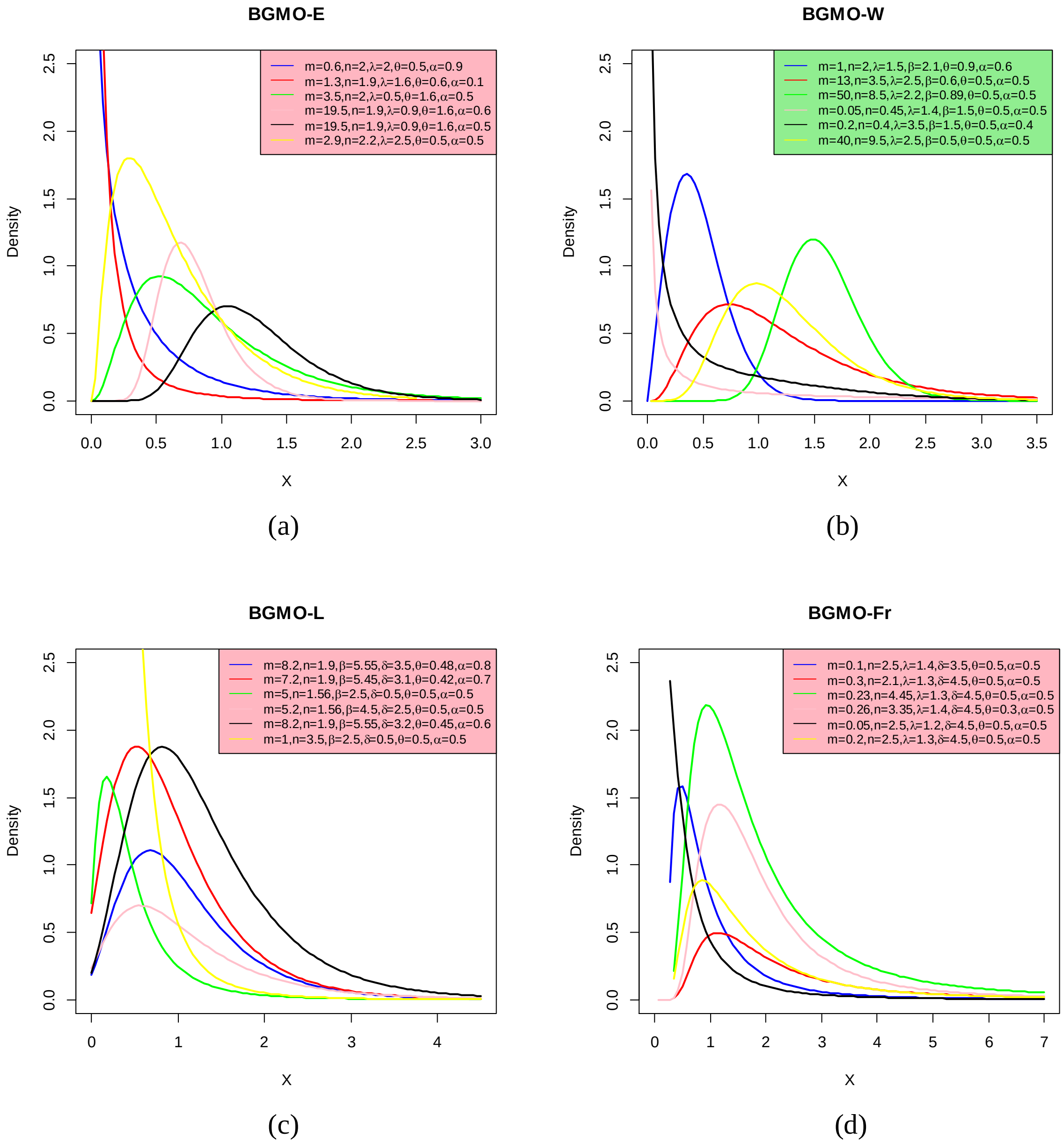

**Fig 1:** Density plots **(a)** $BGMO-E$ , **(b)** $BGMO-W$ , **(c)** $BGMO-L$ and **(d)** $BGMO-Fr$ Distributions:

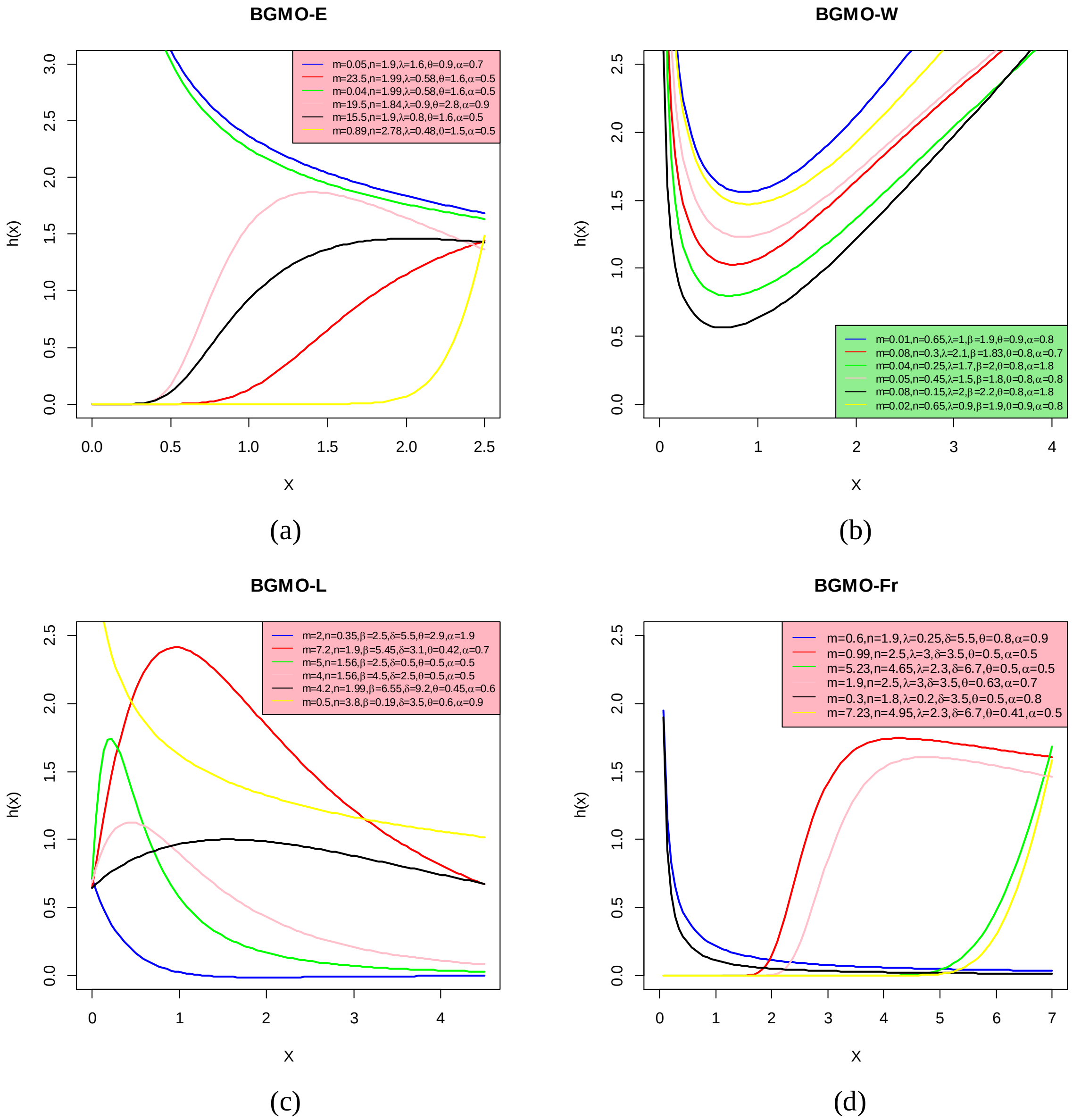

**Fig 2:** Hazard plots plots **(a)** $BGMO-E$, **(b)** $BGMO-W$, **(c)** $BGMO-L$ and **(d)** $BGMO-Fr$ Distributions:

From the plots in figure 1 and 2 it can be seen that the family is very flexible and can offer many different types of shapes of density and hazard rate function including the bath tub shaped free hazard.

## 3. Some special $BGMO-G$ distributions

Some special cases of the $BGMO-G$ family of distributions are presented in this section.

### 3.1 The $BGMO-$ exponential ($BGMO-E$) distribution

Let the base line distribution be exponential with parameter $\lambda > 0$, $g(t:\lambda) = \lambda e^{-\lambda t}, \quad t > 0$ and $G(t:\lambda) = 1 - e^{-\lambda t}, \quad t > 0$ then for the $BGMO-E$ model we get the pdf and cdf respectively as:

$$f^{BGMOE}(t)=\frac{1}{B(m,n)}\frac{\theta\alpha^{\theta}\lambda e^{-\lambda t}[e^{-\lambda t}]^{\theta-1}}{[1-\overline{\alpha}e^{-\lambda t}]^{\theta+1}}\left[1-\left[\frac{\alpha e^{-\lambda t}}{1-\overline{\alpha}e^{-\lambda t}}\right]^{\theta}\right]^{m-1}\left[\left[\frac{\alpha e^{-\lambda t}}{1-\overline{\alpha}e^{-\lambda t}}\right]^{\theta}\right]^{n-1}$$

and $$F^{BGMOE}(t)=I_{1-\left[\frac{\alpha e^{-\lambda t}}{1-\overline{\alpha}e^{-\lambda t}}\right]^{\theta}}(m,n)$$

sf: $$\overline{F}^{BGMOE}(t)=1-I_{1-\left[\frac{\alpha e^{-\lambda t}}{1-\overline{\alpha}e^{-\lambda t}}\right]^{\theta}}(m,n)$$

hrf :

$$h^{BGMOE}(t)=\frac{1}{B(m,n)}\frac{\theta\alpha^{\theta}\lambda e^{-\lambda t}[e^{-\lambda t}]^{\theta-1}}{[1-\overline{\alpha}e^{-\lambda t}]^{\theta+1}\,[1-I_{1-\left[\frac{\alpha e^{-\lambda t}}{1-\overline{\alpha}e^{-\lambda t}}\right]^{\theta}}(m,n)]}$$
$$\times\left[1-\left[\frac{\alpha e^{-\lambda t}}{1-\overline{\alpha}e^{-\lambda t}}\right]^{\theta}\right]^{m-1}\left[\left[\frac{\alpha e^{-\lambda t}}{1-\overline{\alpha}e^{-\lambda t}}\right]^{\theta}\right]^{n-1}$$

rhrf :

$$r^{BGMOE}(t)=\frac{1}{B(m,n)}\frac{\theta\alpha^{\theta}\lambda e^{-\lambda t}[e^{-\lambda t}]^{\theta-1}}{[1-\overline{\alpha}e^{-\lambda t}]^{\theta+1}\,[I_{1-\left[\frac{\alpha e^{-\lambda t}}{1-\overline{\alpha}e^{-\lambda t}}\right]^{\theta}}(m,n)]}$$
$$\times\left[1-\left[\frac{\alpha e^{-\lambda t}}{1-\overline{\alpha}e^{-\lambda t}}\right]^{\theta}\right]^{m-1}\left[\left[\frac{\alpha e^{-\lambda t}}{1-\overline{\alpha}e^{-\lambda t}}\right]^{\theta}\right]^{n-1}$$

chrf: $$H^{BGMOE}(t)=-\log\left[\,1-I_{1-\left[\frac{\alpha e^{-\lambda t}}{1-\overline{\alpha}e^{-\lambda t}}\right]^{\theta}}(m,n)\right]$$

### 3.2 The $BGMO-$Lomax ($BGMO-L$) distribution

Considering the Lomax distribution (Ghitany *et al.* 2007) with pdf and cdf given by $g(t:\beta,\delta)=(\beta/\delta)[1+(t/\delta)]^{-(\beta+1)},\ \ t>0$, and $G(t:\beta,\delta)=1-[1+(t/\delta)]^{-\beta},\ \ \beta>0$ and $\delta>0$ the pdf and cdf of the $BGMO-L$ distribution are given by

$$f^{BGMOL}(t)=\frac{1}{B(m,n)}\frac{\theta\alpha^{\theta}(\beta/\delta)[1+(t/\delta)]^{-(\beta+1)}[[1+(t/\delta)]^{-\beta}]^{\theta-1}}{[1-\overline{\alpha}[1+(t/\delta)]^{-\beta}]^{\theta+1}}$$
$$\times\left[1-\left[\frac{\alpha[1+(t/\delta)]^{-\beta}}{1-\overline{\alpha}[1+(t/\delta)]^{-\beta}}\right]^{\theta}\right]^{m-1}\left[\left[\frac{\alpha[1+(t/\delta)]^{-\beta}}{1-\overline{\alpha}[1+(t/\delta)]^{-\beta}}\right]^{\theta}\right]^{n-1}$$

and $$F^{BGMOL}(t)=I_{1-\left[\frac{\alpha[1+(t/\delta)]^{-\beta}}{1-\overline{\alpha}[1+(t/\delta)]^{-\beta}}\right]^{\theta}}(m,n)$$

sf: $$\overline{F}^{BGMOL}(t)=1-I_{1-\left[\frac{\alpha[1+(t/\delta)]^{-\beta}}{1-\overline{\alpha}[1+(t/\delta)]^{-\beta}}\right]^{\theta}}(m,n)$$

hrf : $h^{BGMOL}(t)$

$$=\frac{1}{B(m,n)}\frac{\theta\alpha^{\theta}(\beta/\delta)[1+(t/\delta)]^{-(\beta+1)}[1+(t/\delta)]^{-\beta(\theta-1)}}{[1-\bar{\alpha}\{1+(t/\delta)\}^{-\beta}]^{\theta+1}}\frac{1}{1-I_{1-\left[\frac{\alpha[1+(t/\delta)]^{-\beta}}{1-\bar{\alpha}[1+(t/\delta)]^{-\beta}}\right]^{\theta}}(m,n)}$$

$$\times\left[1-\left[\frac{\alpha[1+(t/\delta)]^{-\beta}}{1-\bar{\alpha}[1+(t/\delta)]^{-\beta}}\right]^{\theta}\right]^{m-1}\left[\left[\frac{\alpha[1+(t/\delta)]^{-\beta}}{1-\bar{\alpha}[1+(t/\delta)]^{-\beta}}\right]^{\theta}\right]^{n-1}$$

rhrf: $r^{BGMOL}(t)$

$$=\frac{1}{B(m,n)}\frac{\theta\alpha^{\theta}(\beta/\delta)[1+(t/\delta)]^{-(\beta+1)}[1+(t/\delta)]^{-\beta(\theta-1)}}{[1-\bar{\alpha}\{1+(t/\delta)\}^{-\beta}]^{\theta+1}}\frac{1}{I_{1-\left[\frac{\alpha[1+(t/\delta)]^{-\beta}}{1-\bar{\alpha}[1+(t/\delta)]^{-\beta}}\right]^{\theta}}(m,n)}$$

$$\times\left[1-\left[\frac{\alpha[1+(t/\delta)]^{-\beta}}{1-\bar{\alpha}[1+(t/\delta)]^{-\beta}}\right]^{\theta}\right]^{m-1}\left[\left[\frac{\alpha[1+(t/\delta)]^{-\beta}}{1-\bar{\alpha}[1+(t/\delta)]^{-\beta}}\right]^{\theta}\right]^{n-1}$$

chrf: $H^{BGMOL}(t)=-\log\left[1-I_{1-\left[\frac{\alpha[1+(t/\delta)]^{-\beta}}{1-\bar{\alpha}[1+(t/\delta)]^{-\beta}}\right]^{\theta}}(m,n)\right]$

### 3.3 The $BGMO-$Weibull ($BGMO-W$) distribution

Considering the Weibull distribution (Ghitany *et al*. 2005, Zhang and Xie 2007) with parameters $\lambda>0$ and $\beta>0$ having pdf and cdf $g(t)=\lambda\beta t^{\beta-1}e^{-\lambda t^{\beta}}$ and $G(t)=1-e^{-\lambda t^{\beta}}$ respectively we get the pdf and cdf of $BGMO-W$ distribution as

$$f^{BGMOW}(t)=\frac{1}{B(m,n)}\frac{\theta\alpha^{\theta}\lambda\beta t^{\beta-1}e^{-\lambda t^{\beta}}[e^{-\lambda t^{\beta}}]^{\theta-1}}{[1-\bar{\alpha}e^{-\lambda t^{\beta}}]^{\theta+1}}\left[1-\left[\frac{\alpha e^{-\lambda t^{\beta}}}{1-\bar{\alpha}e^{-\lambda t^{\beta}}}\right]^{\theta}\right]^{m-1}\left[\left[\frac{\alpha e^{-\lambda t^{\beta}}}{1-\bar{\alpha}e^{-\lambda t^{\beta}}}\right]^{\theta}\right]^{n-1}$$

and $$F^{BGMOW}(t)=I_{1-\left[\frac{\alpha e^{-\lambda t^{\beta}}}{1-\bar{\alpha}e^{-\lambda t^{\beta}}}\right]^{\theta}}(m,n)$$

sf: $$\bar{F}^{BGMOW}(t)=1-I_{1-\left[\frac{\alpha e^{-\lambda t^{\beta}}}{1-\bar{\alpha}e^{-\lambda t^{\beta}}}\right]^{\theta}}(m,n)$$

hrf: $h^{BGMOW}(t)$

$$=\frac{1}{B(m,n)}\frac{\theta\alpha^{\theta}\lambda\beta t^{\beta-1}e^{-\lambda t^{\beta}}[e^{-\lambda t^{\beta}}]^{\theta-1}}{[1-\bar{\alpha}e^{-\lambda t^{\beta}}]^{\theta+1}}\frac{1}{1-I_{1-\left[\frac{\alpha e^{-\lambda t^{\beta}}}{1-\bar{\alpha}e^{-\lambda t^{\beta}}}\right]^{\theta}}(m,n)}$$

$$\times\left[1-\left[\frac{\alpha e^{-\lambda t^{\beta}}}{1-\bar{\alpha}e^{-\lambda t^{\beta}}}\right]^{\theta}\right]^{m-1}\left[\left[\frac{\alpha e^{-\lambda t^{\beta}}}{1-\bar{\alpha}e^{-\lambda t^{\beta}}}\right]^{\theta}\right]^{n-1}$$

rhrf:

$$r^{BGMOW}(t)=\frac{1}{B(m,n)}\frac{\theta\alpha^{\theta}\lambda\beta\, t^{\beta-1}e^{-\lambda t^{\beta}}[e^{-\lambda t^{\beta}}]^{\theta-1}}{[1-\bar{\alpha}\, e^{-\lambda t^{\beta}}]^{\theta+1}}\frac{1}{I_{1-\left[\frac{\alpha e^{-\lambda t^{\beta}}}{1-\bar{\alpha}\, e^{-\lambda t^{\beta}}}\right]^{\theta}}(m,n)}$$

$$\times\left[1-\left[\frac{\alpha e^{-\lambda t^{\beta}}}{1-\bar{\alpha}\, e^{-\lambda t^{\beta}}}\right]^{\theta}\right]^{m-1}\left[\left[\frac{\alpha e^{-\lambda t^{\beta}}}{1-\bar{\alpha}\, e^{-\lambda t^{\beta}}}\right]^{\theta}\right]^{n-1}$$

chrf: $H^{BGMOW}(t)=-\log\left[\,1-I_{1-\left[\frac{\alpha e^{-\lambda t^{\beta}}}{1-\bar{\alpha}\, e^{-\lambda t^{\beta}}}\right]^{\theta}}(m,n)\right]$

### 3.4 The $BGMO-$ Frechet ( $BGMO-Fr$ ) distribution

Suppose the base line distribution is the Frechet distribution (Krishna *et al.*, 2013) with pdf and cdf given by $g(t)=\lambda\delta^{\lambda}t^{-(\lambda+1)}e^{-(\delta/t)^{\lambda}}$ and $G(t)=e^{-(\delta/t)^{\lambda}}, t>0$ respectively, and then the corresponding pdf and cdf of $BGMO-Fr$ distribution becomes

$$f^{BGMOFr}(t)=\frac{1}{B(m,n)}\frac{\theta\alpha^{\theta}\lambda\delta^{\lambda}t^{-(\lambda+1)}e^{-(\delta/t)^{\lambda}}[1-e^{-(\delta/t)^{\lambda}}]^{\theta-1}}{[1-\bar{\alpha}\,[1-e^{-(\delta/t)^{\lambda}}]]^{\theta+1}}$$

$$\times\left[1-\left[\frac{\alpha\,[1-e^{-(\delta/t)^{\lambda}}]}{1-\bar{\alpha}[1-e^{-(\delta/t)^{\lambda}}]}\right]^{\theta}\right]^{m-1}\left[\left[\frac{\alpha[1-e^{-(\delta/t)^{\lambda}}]}{1-\bar{\alpha}[1-e^{-(\delta/t)^{\lambda}}]}\right]^{\theta}\right]^{n-1}$$

and $F^{BGMOFr}(t)=I_{1-\left[\frac{\alpha\,[1-e^{-(\delta/t)^{\lambda}}]}{1-\bar{\alpha}[1-e^{-(\delta/t)^{\lambda}}]}\right]^{\theta}}(m,n)$

sf: $\bar{F}^{BGMOFr}(t)=1-I_{1-\left[\frac{\alpha\,[1-e^{-(\delta/t)^{\lambda}}]}{1-\bar{\alpha}[1-e^{-(\delta/t)^{\lambda}}]}\right]^{\theta}}(m,n)$

hrf : $h^{BGMOFr}(t)=$

$$\frac{1}{B(m,n)}\frac{\theta\alpha^{\theta}\lambda\delta^{\lambda}t^{-(\lambda+1)}e^{-(\delta/t)^{\lambda}}[1-e^{-(\delta/t)^{\lambda}}]^{\theta-1}}{[1-\bar{\alpha}\,[1-e^{-(\delta/t)^{\lambda}}]]^{\theta+1}}\frac{1}{1-I_{1-\left[\frac{\alpha\,[1-e^{-(\delta/t)^{\lambda}}]}{1-\bar{\alpha}[1-e^{-(\delta/t)^{\lambda}}]}\right]^{\theta}}(m,n)}$$

$$\times\left[1-\left[\frac{\alpha\,[1-e^{-(\delta/t)^{\lambda}}]}{1-\bar{\alpha}[1-e^{-(\delta/t)^{\lambda}}]}\right]^{\theta}\right]^{m-1}\left[\left[\frac{\alpha[1-e^{-(\delta/t)^{\lambda}}]}{1-\bar{\alpha}[1-e^{-(\delta/t)^{\lambda}}]}\right]^{\theta}\right]^{n-1}$$

rhrf : $r^{BGMOFr}(t)=$

$$\frac{1}{B(m,n)}\frac{\theta\alpha^{\theta}\lambda\delta^{\lambda}t^{-(\lambda+1)}e^{-(\delta/t)^{\lambda}}[1-e^{-(\delta/t)^{\lambda}}]^{\theta-1}}{[1-\overline{\alpha}[1-e^{-(\delta/t)^{\lambda}}]]^{\theta+1}}\frac{1}{I_{1-\left[\frac{\alpha[1-e^{-(\delta/t)^{\lambda}}]}{1-\overline{\alpha}[1-e^{-(\delta/t)^{\lambda}}]}\right]^{\theta}}(m,n)}$$

$$\times\left[1-\left[\frac{\alpha[1-e^{-(\delta/t)^{\lambda}}]}{1-\overline{\alpha}[1-e^{-(\delta/t)^{\lambda}}]}\right]^{\theta}\right]^{m-1}\left[\left[\frac{\alpha[1-e^{-(\delta/t)^{\lambda}}]}{1-\overline{\alpha}[1-e^{-(\delta/t)^{\lambda}}]}\right]^{\theta}\right]^{n-1}$$

chrf: $H^{BGMOFr}(t)=-\log[1-I_{1-\left[\frac{\alpha[1-e^{-(\delta/t)^{\lambda}}]}{1-\overline{\alpha}[1-e^{-(\delta/t)^{\lambda}}]}\right]^{\theta}}(m,n)]$

### 3.5 The $BGMO-$Gompertz ($BGMO-Go$) distribution

Next by taking the Gompertz distribution (Gieser et al. 1998) with pdf and cdf $g(t)=\beta e^{\lambda t}e^{-\frac{\beta}{\lambda}(e^{\lambda t}-1)}$ and $G(t)=1-e^{-\frac{\beta}{\lambda}(e^{\lambda t}-1)}$, $\beta>0,\ \lambda>0,\ t>0$ respectively, we get the pdf and cdf of $BGMO-Go$ distribution as

$f^{BGMOGo}(t)=$

$$\frac{1}{B(m,n)}\frac{\theta\alpha^{\theta}\beta e^{\lambda t}e^{-\frac{\beta}{\lambda}(e^{\lambda t}-1)}[e^{-\frac{\beta}{\lambda}(e^{\lambda t}-1)}]^{\theta-1}}{[1-\overline{\alpha}e^{-\frac{\beta}{\lambda}(e^{\lambda t}-1)}]^{\theta+1}}\left[1-\left[\frac{\alpha e^{-\frac{\beta}{\lambda}(e^{\lambda t}-1)}}{1-\overline{\alpha}e^{-\frac{\beta}{\lambda}(e^{\lambda t}-1)}}\right]^{\theta}\right]^{m-1}\left[\left[\frac{\alpha e^{-\frac{\beta}{\lambda}(e^{\lambda t}-1)}}{1-\overline{\alpha}e^{-\frac{\beta}{\lambda}(e^{\lambda t}-1)}}\right]^{\theta}\right]^{n-1}$$

and $F^{BGMOGo}(t)=I_{1-\left[\frac{\alpha e^{-\frac{\beta}{\lambda}(e^{\lambda t}-1)}}{1-\overline{\alpha}e^{-\frac{\beta}{\lambda}(e^{\lambda t}-1)}}\right]^{\theta}}(m,n)$

sf: $\overline{F}^{BGMOGo}(t)=1-I_{1-\left[\frac{\alpha e^{-\frac{\beta}{\lambda}(e^{\lambda t}-1)}}{1-\overline{\alpha}e^{-\frac{\beta}{\lambda}(e^{\lambda t}-1)}}\right]^{\theta}}(m,n)$

hrf: $h^{BGMOGo}(t)=$

$$\frac{1}{B(m,n)}\frac{\theta\alpha^{\theta}\beta e^{\lambda t}e^{-\frac{\beta}{\lambda}(e^{\lambda t}-1)}[e^{-\frac{\beta}{\lambda}(e^{\lambda t}-1)}]^{\theta-1}}{[1-\overline{\alpha}e^{-\frac{\beta}{\lambda}(e^{\lambda t}-1)}]^{\theta+1}}\frac{1}{1-I_{1-\left[\frac{\alpha e^{-\frac{\beta}{\lambda}(e^{\lambda t}-1)}}{1-\overline{\alpha}e^{-\frac{\beta}{\lambda}(e^{\lambda t}-1)}}\right]^{\theta}}(m,n)}$$

$$\times\left[1-\left[\frac{\alpha e^{-\frac{\beta}{\lambda}(e^{\lambda t}-1)}}{1-\overline{\alpha}e^{-\frac{\beta}{\lambda}(e^{\lambda t}-1)}}\right]^{\theta}\right]^{m-1}\left[\left[\frac{\alpha e^{-\frac{\beta}{\lambda}(e^{\lambda t}-1)}}{1-\overline{\alpha}e^{-\frac{\beta}{\lambda}(e^{\lambda t}-1)}}\right]^{\theta}\right]^{n-1}$$

rhrf: $r^{BGMOGo}(t)=$

$$\frac{1}{B(m,n)}\frac{\theta\alpha^{\theta}\beta e^{\lambda t}e^{-\frac{\beta}{\lambda}(e^{\lambda t}-1)}[e^{-\frac{\beta}{\lambda}(e^{\lambda t}-1)}]^{\theta-1}}{[1-\overline{\alpha}\,e^{-\frac{\beta}{\lambda}(e^{\lambda t}-1)}]^{\theta+1}}\frac{1}{I_{1-\left[\frac{\alpha e^{-\frac{\beta}{\lambda}(e^{\lambda t}-1)}}{1-\overline{\alpha}\,e^{-\frac{\beta}{\lambda}(e^{\lambda t}-1)}}\right]^{\theta}}(m,n)}$$

$$\times\left[1-\left[\frac{\alpha e^{-\frac{\beta}{\lambda}(e^{\lambda t}-1)}}{1-\overline{\alpha}\,e^{-\frac{\beta}{\lambda}(e^{\lambda t}-1)}}\right]^{\theta}\right]^{m-1}\left[\left[\frac{\alpha e^{-\frac{\beta}{\lambda}(e^{\lambda t}-1)}}{1-\overline{\alpha}\,e^{-\frac{\beta}{\lambda}(e^{\lambda t}-1)}}\right]^{\theta}\right]^{n-1}$$

chrf: $$H^{BGMOGo}(t)=-\log\,[\,1-I_{1-\left[\frac{\alpha e^{-\frac{\beta}{\lambda}(e^{\lambda t}-1)}}{1-\overline{\alpha}\,e^{-\frac{\beta}{\lambda}(e^{\lambda t}-1)}}\right]^{\theta}}(m,n)]$$

### 3.6 The $BGMO-$ Extended Weibull ( $BGMO-EW$ ) distribution

The pdf and the cdf of the extended Weibull (*EW*) distributions of Gurvich *et al.* (1997) is given by $g(t:\delta,\vartheta)=\delta\exp[-\delta Z(t:\vartheta)]z(t:\vartheta)$ and $G(t:\delta,\vartheta)=1-\exp[-\delta Z(t:\vartheta)]$, $t\in D\subseteq R_{+},\delta>0$ where $Z(t:\vartheta)$ is a non-negative monotonically increasing function which depends on the parameter vector $\vartheta$. and $z(t:\vartheta)$ is the derivative of $Z(t:\vartheta)$.

By considering *EW* as the base line distribution we derive pdf and cdf of the $BGMO-EW$ as

$$f^{BGMOEW}(t)=\frac{1}{B(m,n)}\frac{\theta\alpha^{\theta}\delta\exp[-\delta Z(t:\vartheta)]z(t:\vartheta)[\exp\{-\delta Z(t:\vartheta)\}]^{\theta-1}}{[1-\overline{\alpha}\exp\{-\delta Z(t:\vartheta)\}]^{\theta+1}}$$

$$\times\left[1-\left[\frac{\alpha\exp[-\delta Z(t:\vartheta)]}{1-\overline{\alpha}\exp[-\delta Z(t:\vartheta)]}\right]^{\theta}\right]^{m-1}\left[\left[\frac{\alpha\exp[-\delta Z(t:\vartheta)]}{1-\overline{\alpha}\exp[-\delta Z(t:\vartheta)]}\right]^{\theta}\right]^{n-1}$$

and $$F^{BGMOEW}(t)=I_{1-\left[\frac{\alpha\exp[-\delta Z(t:\vartheta)]}{1-\overline{\alpha}\exp[-\delta Z(t:\vartheta)]}\right]^{\theta}}(m,n)$$

Important models can be seen as particular cases with different choices of $Z(t:\vartheta)$:

(i) $Z(t:\vartheta)=t$: exponential distribution.

(ii) $Z(t:\vartheta)=t^{2}$: Rayleigh (Burr type-X) distribution.

(iii) $Z(t:\vartheta)=\log(t/k)$: Pareto distribution

(iv) $Z(t:\vartheta)=\beta^{-1}[\exp(\beta t)-1]$: Gompertz distribution.

sf: $$\overline{F}^{BGMOEW}(t)=1-I_{1-\left[\frac{\alpha\exp[-\delta Z(t:\vartheta)]}{1-\overline{\alpha}\exp[-\delta Z(t:\vartheta)]}\right]^{\theta}}(m,n)$$

hrf: $h^{BGMOEW}(t) =$

$$\frac{1}{B(m,n)} \frac{\theta\alpha^{\theta}\delta\exp[-\delta Z(t:\vartheta)]z(t:\vartheta)[\exp\{-\delta Z(t:\vartheta)\}]^{\theta-1}}{[1-\bar{\alpha}\exp\{-\delta Z(t:\vartheta)\}]^{\theta+1}} \frac{1}{1-I_{1-\left[\frac{\alpha\exp[-\delta Z(t:\vartheta)]}{1-\bar{\alpha}\exp[-\delta Z(t:\vartheta)]}\right]^{\theta}}(m,n)}$$

$$\times\left[1-\left[\frac{\alpha\exp[-\delta Z(t:\vartheta)]}{1-\bar{\alpha}\exp[-\delta Z(t:\vartheta)]}\right]^{\theta}\right]^{m-1}\left[\left[\frac{\alpha\exp[-\delta Z(t:\vartheta)]}{1-\bar{\alpha}\exp[-\delta Z(t:\vartheta)]}\right]^{\theta}\right]^{n-1}$$

rhrf: $r^{BGMOEW}(t) =$

$$\frac{1}{B(m,n)} \frac{\theta\alpha^{\theta}\delta\exp[-\delta Z(t:\vartheta)]z(t:\vartheta)[\exp\{-\delta Z(t:\vartheta)\}]^{\theta-1}}{[1-\bar{\alpha}\exp\{-\delta Z(t:\vartheta)\}]^{\theta+1}} \frac{1}{1-I_{1-\left[\frac{\alpha\exp[-\delta Z(t:\vartheta)]}{1-\bar{\alpha}\exp[-\delta Z(t:\vartheta)]}\right]^{\theta}}(m,n)}$$

$$\times\left[1-\left[\frac{\alpha\exp[-\delta Z(t:\vartheta)]}{1-\bar{\alpha}\exp[-\delta Z(t:\vartheta)]}\right]^{\theta}\right]^{m-1}\left[\left[\frac{\alpha\exp[-\delta Z(t:\vartheta)]}{1-\bar{\alpha}\exp[-\delta Z(t:\vartheta)]}\right]^{\theta}\right]^{n-1}$$

chrf: $$H^{BGMOEW}(t) = -\log\left[1-I_{1-\left[\frac{\alpha\exp[-\delta Z(t:\vartheta)]}{1-\bar{\alpha}\exp[-\delta Z(t:\vartheta)]}\right]^{\theta}}(m,n)\right]$$

### 3.7 The $BGMO-$ Extended Modified Weibull ( $BGMO-EMW$ ) distribution

The modified Weibull (MW) distribution (Sarhan and Zaindin 2013) with cdf and pdf is given by

$G(t;\sigma,\beta,\gamma) = 1-\exp[-\sigma t-\beta t^{\gamma}]$ , $t>0, \gamma>0, \sigma,\beta\geq 0, \sigma+\beta>0$ and

$g(t;\sigma,\beta,\gamma) = (\sigma+\beta\gamma t^{\gamma-1})\exp[-\sigma t-\beta t^{\gamma}]$ respectively.

The corresponding pdf and cdf of $BGMO-EMW$ are given by

$$f^{BGMOEMW}(t) = \frac{1}{B(m,n)} \frac{\theta\alpha^{\theta}(\sigma+\beta\gamma t^{\gamma-1})\exp[-\sigma t-\beta t^{\gamma}][\exp\{-\sigma t-\beta t^{\gamma}\}]^{\theta-1}}{[1-\bar{\alpha}\exp\{-\sigma t-\beta t^{\gamma}\}]]^{\theta+1}}$$

$$\times\left[1-\left[\frac{\alpha\exp[-\sigma t-\beta t^{\gamma}]}{1-\bar{\alpha}\exp[-\sigma t-\beta t^{\gamma}]}\right]^{\theta}\right]^{m-1}\left[\left[\frac{\alpha\exp[-\sigma t-\beta t^{\gamma}]}{1-\bar{\alpha}\exp[-\sigma t-\beta t^{\gamma}]}\right]^{\theta}\right]^{n-1}$$

and $$F^{BGMOEMW}(t) = I_{1-\left[\frac{\alpha\exp[-\sigma t-\beta t^{\gamma}]}{1-\bar{\alpha}\exp[-\sigma t-\beta t^{\gamma}]}\right]^{\theta}}(m,n)$$

sf: $$\bar{F}^{BGMOEMW}(t) = 1-I_{1-\left[\frac{\alpha\exp[-\sigma t-\beta t^{\gamma}]}{1-\bar{\alpha}\exp[-\sigma t-\beta t^{\gamma}]}\right]^{\theta}}(m,n)$$

hrf : $h^{BGMOEMW}(t)=$

$$\frac{1}{B(m,n)}\frac{\theta\alpha^{\theta}(\sigma+\beta\gamma t^{\gamma-1})\exp[-\sigma t-\beta t^{\gamma}][\exp\{-\sigma t-\beta t^{\gamma}\}]^{\theta-1}}{[1-\bar{\alpha}\exp\{-\sigma t-\beta t^{\gamma}\}]]^{\theta+1}}\frac{1}{1-I_{1-\left[\frac{\alpha\exp[-\sigma t-\beta t^{\gamma}]}{1-\bar{\alpha}\exp[-\sigma t-\beta t^{\gamma}]}\right]^{\theta}}(m,n)}$$

$$\times\left[1-\left[\frac{\alpha\exp[-\sigma t-\beta t^{\gamma}]}{1-\bar{\alpha}\exp[-\sigma t-\beta t^{\gamma}]}\right]^{\theta}\right]^{m-1}\left[\left[\frac{\alpha\exp[-\sigma t-\beta t^{\gamma}]}{1-\bar{\alpha}\exp[-\sigma t-\beta t^{\gamma}]}\right]^{\theta}\right]^{n-1}$$

rhrf: $r^{BGMOEMW}(t)=$

$$\frac{1}{B(m,n)}\frac{\theta\alpha^{\theta}(\sigma+\beta\gamma t^{\gamma-1})\exp[-\sigma t-\beta t^{\gamma}][\exp\{-\sigma t-\beta t^{\gamma}\}]^{\theta-1}}{[1-\bar{\alpha}\exp\{-\sigma t-\beta t^{\gamma}\}]]^{\theta+1}}\frac{1}{I_{1-\left[\frac{\alpha\exp[-\sigma t-\beta t^{\gamma}]}{1-\bar{\alpha}\exp[-\sigma t-\beta t^{\gamma}]}\right]^{\theta}}(m,n)}$$

$$\times\left[1-\left[\frac{\alpha\exp[-\sigma t-\beta t^{\gamma}]}{1-\bar{\alpha}\exp[-\sigma t-\beta t^{\gamma}]}\right]^{\theta}\right]^{m-1}\left[\left[\frac{\alpha\exp[-\sigma t-\beta t^{\gamma}]}{1-\bar{\alpha}\exp[-\sigma t-\beta t^{\gamma}]}\right]^{\theta}\right]^{n-1}$$

chrf: $$H^{BGMOEMW}(t)=-\log\left[1-I_{1-\left[\frac{\alpha\exp[-\sigma t-\beta t^{\gamma}]}{1-\bar{\alpha}\exp[-\sigma t-\beta t^{\gamma}]}\right]^{\theta}}(m,n)\right]$$

### 3.8 The $BGMO-$ Extended Exponentiated Pareto ( $BGMO-EEP$ ) distribution

The pdf and cdf of the exponentiated Pareto distribution, of Nadarajah (2005), are given respectively by $g(t)=\gamma k\theta^{k}t^{-(k+1)}[1-(\theta/t)^{k}]^{\gamma-1}$ and $G(t)=[1-(\theta/t)^{k}]^{\gamma}$, $x>\theta$ and $\theta,k,\gamma>0$. Thus the pdf and the cdf of $BGMO-EEP$ distribution are given by

$$f^{BGMOEEP}(t)=\frac{1}{B(m,n)}\frac{\theta\alpha^{\theta}\gamma k\theta^{k}t^{-(k+1)}[1-(\theta/t)^{k}]^{\gamma-1}[1-\{1-(\theta/t)^{k}\}^{\gamma}]^{\theta-1}}{[1-\bar{\alpha}[1-\{1-(\theta/t)^{k}\}^{\gamma}]]^{\theta+1}}$$

$$\times\left[1-\left[\frac{\alpha[1-\{1-(\theta/t)^{k}\}^{\gamma}]}{1-\bar{\alpha}[1-\{1-(\theta/t)^{k}\}^{\gamma}]}\right]^{\theta}\right]^{m-1}\left[\left[\frac{\alpha[1-\{1-(\theta/t)^{k}\}^{\gamma}]}{1-\bar{\alpha}[1-\{1-(\theta/t)^{k}\}^{\gamma}]}\right]^{\theta}\right]^{n-1}$$

and $$F^{BGMOEEP}(t)=I_{1-\left[\frac{\alpha[1-\{1-(\theta/t)^{k}\}^{\gamma}]}{1-\bar{\alpha}[1-\{1-(\theta/t)^{k}\}^{\gamma}]}\right]^{\theta}}(m,n)$$

sf: $$\bar{F}^{BGMOEEP}(t)=1-I_{1-\left[\frac{\alpha[1-\{1-(\theta/t)^{k}\}^{\gamma}]}{1-\bar{\alpha}[1-\{1-(\theta/t)^{k}\}^{\gamma}]}\right]^{\theta}}(m,n)$$

hrf: $h^{BGMOEEP}(t)=$

$$\frac{1}{B(m,n)}\frac{\theta\alpha^{\theta}\gamma k\theta^{k}t^{-(k+1)}[1-(\theta/t)^{k}]^{\gamma-1}[1-\{1-(\theta/t)^{k}\}^{\gamma}]^{\theta-1}}{[1-\bar{\alpha}[1-\{1-(\theta/t)^{k}\}^{\gamma}]]^{\theta+1}}\frac{1}{1-I_{1-\left[\frac{\alpha[1-\{1-(\theta/t)^{k}\}^{\gamma}]}{1-\bar{\alpha}[1-\{1-(\theta/t)^{k}\}^{\gamma}]}\right]^{\theta}}(m,n)}$$

$$\times\left[1-\left[\frac{\alpha[1-\{1-(\theta/t)^{k}\}^{\gamma}]}{1-\bar{\alpha}[1-\{1-(\theta/t)^{k}\}^{\gamma}]}\right]^{\theta}\right]^{m-1}\left[\left[\frac{\alpha[1-\{1-(\theta/t)^{k}\}^{\gamma}]}{1-\bar{\alpha}[1-\{1-(\theta/t)^{k}\}^{\gamma}]}\right]^{\theta}\right]^{n-1}$$

rhrf: $r^{BGMOEEP}(t) =$

$$\frac{1}{B(m,n)}\frac{\theta\alpha^{\theta}\gamma k\theta^{k}t^{-(k+1)}[1-(\theta/t)^{k}]^{\gamma-1}[1-\{1-(\theta/t)^{k}\}^{\gamma}]^{\theta-1}}{[1-\bar{\alpha}[1-\{1-(\theta/t)^{k}\}^{\gamma}]]^{\theta+1}}\frac{1}{I_{1-\left[\frac{\alpha[1-\{1-(\theta/t)^{k}\}^{\gamma}]}{1-\bar{\alpha}[1-\{1-(\theta/t)^{k}\}^{\gamma}]}\right]^{\theta}}(m,n)}$$

$$\times\left[1-\left[\frac{\alpha[1-\{1-(\theta/t)^{k}\}^{\gamma}]}{1-\bar{\alpha}[1-\{1-(\theta/t)^{k}\}^{\gamma}]}\right]^{\theta}\right]^{m-1}\left[\left[\frac{\alpha[1-\{1-(\theta/t)^{k}\}^{\gamma}]}{1-\bar{\alpha}[1-\{1-(\theta/t)^{k}\}^{\gamma}]}\right]^{\theta}\right]^{n-1}$$

chrf: $H^{BGMOEEP}(t) = -\log\left[\,1-I_{1-\left[\frac{\alpha[1-\{1-(\theta/t)^{k}\}^{\gamma}]}{1-\bar{\alpha}[1-\{1-(\theta/t)^{k}\}^{\gamma}]}\right]^{\theta}}(m,n)\right]$

## 4. General results for the Beta Generalized Marshall-Olkin ($BGMO-G$) family of distributions

In this section we derive some general results for the proposed $BGMO-G$ family.

### 4.1 Expansions

By using binomial expansion in (6), we obtain

$$f^{BGMOG}(t;\alpha,\theta,m,n)$$

$$=\frac{1}{B(m,n)}\frac{\theta\alpha^{\theta}g(t)\bar{G}(t)^{\theta-1}}{[1-\bar{\alpha}\bar{G}(t)]^{\theta+1}}\left[1-\left[\frac{\alpha\bar{G}(t)}{1-\bar{\alpha}\bar{G}(t)}\right]^{\theta}\right]^{m-1}\left[\left[\frac{\alpha\bar{G}(t)}{1-\bar{\alpha}\bar{G}(t)}\right]^{\theta}\right]^{n-1}$$

$$=\frac{\theta}{B(m,n)}\frac{\alpha g(t)}{[1-\bar{\alpha}\bar{G}(t)]^{2}}\left(\frac{\alpha\bar{G}(t)}{[1-\bar{\alpha}\bar{G}(t)]}\right)^{\theta-1}\left[1-\left[\frac{\alpha\bar{G}(t)}{1-\bar{\alpha}\bar{G}(t)}\right]^{\theta}\right]^{m-1}\left[\left[\frac{\alpha\bar{G}(t)}{1-\bar{\alpha}\bar{G}(t)}\right]^{\theta}\right]^{n-1}$$

$$=\frac{\theta}{B(m,n)}f^{MO}(t,\alpha)[\bar{F}^{MO}(t,\alpha)]^{\theta-1}[1-\{\bar{F}^{MO}(t,\alpha)\}^{\theta}]^{m-1}[\bar{F}^{MO}(t,\alpha)]^{\theta(n-1)}$$

$$=\frac{\theta}{B(m,n)}f^{MO}(t,\alpha)[\bar{F}^{MO}(t,\alpha)]^{\theta n-1}[1-\{\bar{F}^{MO}(t,\alpha)\}^{\theta}]^{(m-1)}$$

$$=\frac{\theta}{B(m,n)}f^{MO}(t,\alpha)[\bar{F}^{MO}(t,\alpha)]^{\theta n-1}\sum_{j=0}^{m-1}\binom{m-1}{j}(-1)^{j}[\bar{F}^{MO}(t;\alpha)]^{\theta j}$$

$$=\frac{\theta}{B(m,n)}f^{MO}(t,\alpha)\sum_{j=0}^{m-1}\binom{m-1}{j}(-1)^{j}[\bar{F}^{MO}(t;\alpha)]^{\theta(j+n)-1}$$

$$=f^{MO}(t,\alpha)\sum_{j=0}^{m-1}\delta_{j}[\bar{F}^{MO}(t;\alpha)]^{\theta(j+n)-1} \qquad (10)$$

$$=\sum_{j=0}^{m-1}\frac{\delta_{j}}{\theta(j+n)}\frac{d}{dt}[\bar{F}^{MO}(t;\alpha)]^{\theta(j+n)}$$

$$= \sum_{j=0}^{m-1} \delta'_j \frac{d}{dt} [\bar{F}^{MO}(t;\alpha)]^{\theta(j+n)} \quad (11)$$

$$= \sum_{j=0}^{m-1} \delta'_j \frac{d}{dt} [\bar{F}^{MO}(t;\alpha(\theta(j+n)))]$$

$$= \sum_{j=0}^{m-1} \delta'_j \, f^{MO}(t;\alpha(\theta(j+n))) \quad (12)$$

Where $\delta'_j = \frac{(-1)^{j+1}}{B(m,n)\,(j+n)} \binom{m-1}{j}$ and $\delta_j = \delta'_j \theta(j+n)$

Alternatively, we can expand the pdf as

$$= f^{MO}(t,\alpha) \sum_{j=0}^{m-1} \delta_j [\bar{F}^{MO}(t;\alpha)]^{\theta(j+n)-1}$$

$$= f^{MO}(t,\alpha) \sum_{j=0}^{m-1} \delta_j \sum_{l=0}^{\theta(j+n)-1} \binom{\theta(j+n)-1}{l} (-1)^l [F^{MO}(t;\alpha)]^l$$

$$= f^{MO}(t,\alpha) \sum_{l=0}^{\theta(j+n)-1} \varphi_l [F^{MO}(t;\alpha)]^l \quad (13)$$

$$= \sum_{l=0}^{\theta(j+n)-1} \frac{\varphi_l}{l+1} \frac{d}{dt} [F^{MO}(t;\alpha)]^{l+1}$$

$$= \sum_{l=0}^{\theta(j+n)-1} \varphi'_l \frac{d}{dt} [F^{MO}(t;\alpha)]^{l+1}$$

Where $\varphi_l = \sum_{j=0}^{\infty} \delta_j (-1)^l \binom{\theta(j+n)-1}{l}$ and $\varphi'_l = \frac{1}{l+1} \sum_{j=0}^{\infty} \delta_j (-1)^l \binom{\theta(j+n)-1}{l}$

We can expand the cdf as (see "Incomplete Beta Function" *From Math World*—A Wolfram Web Resource. http://mathworld. Wolfram.com/ Incomplete Beta Function. html)

$B(z;a,b) = B_z(a,b) = z^a \sum_{n=0}^{\infty} \frac{(1-b)_n}{n!\,(a+n)} z^n$ Where $(x)_n$ is a Pochhammer symbol.

$$= z^a \sum_{n=0}^{\infty} \frac{(-1)^n \, (b-1)!}{n!\,(b-n-1)!\,(a+n)} z^n$$

$$= z^a \sum_{n=0}^{\infty} \binom{b-1}{n} \frac{(-1)^n}{(a+n)} z^n \quad (14)$$

$$F^{BGMOG}(t) = I_{1-\left[\frac{\alpha \bar{G}(t)}{1-\bar{\alpha}\,\bar{G}(t)}\right]^{\theta}}(m,n)$$

Using (14) in (7) we have

$$=\frac{1}{B(m,n)}\left[1-\left(\frac{\alpha\overline{G}(t)}{1-\overline{\alpha}\,\overline{G}(t)}\right)^{\theta}\right]^{m}\sum_{i=0}^{\infty}\binom{n-1}{i}\frac{(-1)^{i}}{(m+i)}\left[1-\left(\frac{\alpha\overline{G}(t)}{1-\overline{\alpha}\,\overline{G}(t)}\right)^{\theta}\right]^{i}$$

$$=\frac{1}{B(m,n)}[1-\{\overline{F}^{MO}(t,\alpha)\}^{\theta}]^{m}\sum_{i=0}^{\infty}\binom{n-1}{i}\frac{(-1)^{i}}{(m+i)}[1-\{\overline{F}^{MO}(t,\alpha)\}^{\theta}]^{i}$$

$$=\sum_{i=0}^{\infty}\frac{(-1)^{i}}{B(m,n)(m+i)}\binom{n-1}{i}[1-\{\overline{F}^{MO}(t,\alpha)\}^{\theta}]^{m+i}$$

$$=\sum_{i=0}^{\infty}\frac{(-1)^{i}}{B(m,n)(m+i)}\binom{n-1}{i}\sum_{j=0}^{\infty}\binom{m+i}{j}(-1)^{j}[\overline{F}^{MO}(t,\alpha)]^{\theta j}$$

$$=\sum_{i,j=0}^{\infty}\frac{(-1)^{i+j}}{B(m,n)(m+i)}\binom{n-1}{i}\binom{m+i}{j}\sum_{k=0}^{j}\binom{\theta j}{k}(-1)^{k}[F^{MO}(t,\alpha)]^{k}$$

$$=\sum_{i,j=0}^{\infty}\sum_{k=0}^{j}\frac{(-1)^{i+j+k}}{B(m,n)(m+i)}\binom{n-1}{i}\binom{m+i}{j}\binom{\theta j}{k}[F^{MO}(t,\alpha)]^{k}$$

By exchanging the indices $j$ and $k$ in the sum symbol, we have

$$F^{BGMO}(t;\alpha,\theta,m,n)=\sum_{i,k=0}^{\infty}\sum_{j=k}^{\infty}\frac{(-1)^{i+j+k}}{B(m,n)(m+i)}\binom{n-1}{i}\binom{m+i}{j}\binom{\theta j}{k}[F^{MO}(t,\alpha)]^{k}$$

and then

$$F^{BGMO}(t;\alpha,\theta,m,n)=\sum_{k=0}^{\infty}\chi_{k}\,F^{MO}(t,\alpha)^{k}\qquad(15)$$

Where $\chi_{k}=\sum_{i=0}^{\infty}\sum_{j=k}^{\infty}\frac{(-1)^{i+j+k}}{B(m,n)(m+i)}\binom{n-1}{i}\binom{m+i}{j}\binom{\theta j}{k}$

Similarly an expansion for the cdf of $BGMO-G$ can be derives as

$$F^{BGMOG}(t;\alpha,\theta,m,n)=I_{1-[\overline{F}^{MO}(t,\alpha)]^{\theta}}(m,n)$$

$$=\sum_{p=m}^{m+n-1}\binom{m+n-1}{p}[1-\{\overline{F}^{MO}(t,\alpha)\}^{\theta}]^{p}[\overline{F}^{MO}(t,\alpha)^{\theta}]^{m+n-1-p}$$

$$=\sum_{p=m}^{m+n-1}\binom{m+n-1}{p}\sum_{q=0}^{p}\binom{p}{q}(-1)^{q}\{\overline{F}^{MO}(t,\alpha)\}^{\theta q}[\overline{F}^{MO}(t,\alpha)^{\theta}]^{m+n-1-p}$$

$$=\sum_{p=m}^{m+n-1}\sum_{q=0}^{p}(-1)^{q}\binom{p}{q}\binom{m+n-1}{p}[\overline{F}^{MO}(t,\alpha)]^{\theta(m+n-1-p+q)}$$

$$=\sum_{p=m}^{m+n-1}\sum_{q=0}^{p}(-1)^{q}\binom{p}{q}\binom{m+n-1}{p}\sum_{r=0}^{\theta(m+n-1-p+q)}\binom{\theta(m+n-1-p+q)}{r}(-1)^{r}[F^{MO}(t,\alpha)]^{r}$$

$$=\sum_{p=m}^{m+n-1}\sum_{q=0}^{p}\sum_{r=0}^{\theta(m+n-1-p+q)}(-1)^{q+r}\binom{p}{q}\binom{m+n-1}{p}\binom{\theta(m+n-1-p+q)}{r}[F^{MO}(t,\alpha)]^{r}$$

$$= \sum_{p=m}^{m+n-1} \sum_{q=0}^{p} \sum_{r=0}^{\theta(m+n-1-p+q)} \psi_{p,q,r} [F^{MO}(t,\alpha)]^{r} \tag{16}$$

Where $\psi_{p,q,r} = (-1)^{q+r} \binom{p}{q} \binom{m+n-1}{p} \binom{\theta(m+n-1-p+q)}{r}$

**4.2 Order statistics**

Suppose $T_1, T_2, ... T_n$ is a random sample from any $BGMO-G$ distribution. Let $T_{r:n}$ denote the $r^{th}$ order statistics. The pdf of $T_{r:n}$ can be expressed as

$$f_{r:n}(t) = \frac{n!}{(r-1)!\,(n-r)!} f^{BGMOG}(t)\; F^{BGMOG}(t)^{r-1} \{1 - F^{BGMOG}(t)\}^{n-r}$$

$$= \frac{n!}{(r-1)!\,(n-r)!} \sum_{j=0}^{n-r} (-1)^{j} \binom{n-r}{j} f^{BGMOG}(t) F^{BGMOG}(t)^{j+r-1}$$

Now using the general expansion of the $BGMO-G$ distribution pdf and cdf we get the pdf of the $r^{th}$ order statistics for of the $BGMO-G$ is given by

$$f_{r:n}(t) = \frac{n!}{(r-1)!\,(n-r)!} \sum_{j=0}^{n-r} (-1)^{j} \binom{n-r}{j} \left\{ f^{MO}(t,\alpha) \sum_{l=0}^{\theta(j+n)-1} \varphi_l [F^{MO}(t;\alpha)]^{l} \right\} \left\{ \sum_{k=0}^{\infty} \chi_k \; F^{MO}(t,\alpha)^{k} \right\}^{j+r-1}$$

Where $\varphi_l$ and $\chi_k$ defined in above

Now $\left\{ \sum_{k=0}^{\infty} \chi_k \; F^{MO}(t,\alpha)^{k} \right\}^{j+r-1} = \sum_{k=0}^{\infty} d_{j+r-1,k} [F^{MO}(t,\alpha)]^{k}$

Where $d_{j+r-1,k} = \frac{1}{k\,\chi_0} \sum_{c=1}^{k} [c(j+r)-k] \chi_c \, d_{j+r-1,\; k-c}$ (Nadarajah et. al 2015)

Therefore the density function of the $r^{th}$ order statistics of $BGMO-G$ distribution can be expressed as

$$f_{r:n}(t) = \frac{n!}{(r-1)!\,(n-r)!} \sum_{j=0}^{n-r} (-1)^{j} \binom{n-r}{j} \left\{ f^{MO}(t,\alpha) \sum_{l=0}^{\theta(j+n)-1} \varphi_l [F^{MO}(t;\alpha)]^{l} \right\} \left( \sum_{k=0}^{\infty} d_{j+r-1,k} [F^{MO}(t,\alpha)]^{k} \right)$$

$$= \frac{n!}{(r-1)!\,(n-r)!} \sum_{j=0}^{n-r} (-1)^{j} \binom{n-r}{j} \left\{ f^{MO}(t,\alpha) \sum_{l=0}^{\theta(j+n)-1} \sum_{k=0}^{\infty} \varphi_l \; d_{j+r-1,k} \; [F^{MO}(t,\alpha)]^{k+l} \right\}$$

$$= f^{MO}(t,\alpha) \sum_{l=0}^{\theta(j+n)-1} \sum_{k=0}^{\infty} \xi_{l,k} [F^{MO}(t,\alpha)]^{k+l} \tag{17}$$

Where $\xi_{l,k} = \dfrac{n!}{(r-1)!\,(n-r)!}\sum_{j=0}^{n-r}(-1)^{j}\binom{n-r}{j}\varphi_l\, d_{j+r-1,k}$ and $\varphi_l$ and $d_{j+r-1,k}$ defined above.

**4.3 Probability weighted moments**

The probability weighted moments (PWMs), first proposed by Greenwood et al. (1979), are expectations of certain functions of a random variable whose mean exists. The $(p,q,r)^{th}$ PWM of $T$ is defined by

$$\Gamma_{p,q,r} = \int_{-\infty}^{\infty} t^p F(t)^q\,[1-F(t)]^r\, f(t)\,dt$$

From equations (10) and (13) the $s^{th}$ moment of $T$ can be written either as

$$E(T^s) = \int_{-\infty}^{\infty} t^s\;f^{BGMOG}(t;\alpha,\theta,m,n)\,dt$$

$$= \sum_{j=0}^{m-1}\delta_j \int_{-\infty}^{\infty} t^s\,[\overline{F}^{MO}(t;\alpha)]^{\,\theta(j+n)-1}\, f^{MO}(t,\alpha)\,dt$$

$$= \sum_{j=0}^{m-1}\delta_j \int_{-\infty}^{\infty} t^s\,[\alpha\,\overline{G}(t)/1-\overline{\alpha}\,\overline{G}(t)]^{\theta(j+n)-1}\,[\alpha\, g(t)/[1-\overline{\alpha}\,\overline{G}(t)]^2]\,dt$$

$$= \sum_{j=0}^{m-1}\delta_j\;\Gamma_{s,0,\,\theta(j+n)-1}$$

or $E(T^s) = \sum_{l=0}^{\theta(j+n)-1}\varphi_l \int_{-\infty}^{\infty} t^s\,[F^{MO}(t;\alpha)]^l\, f^{MO}(t,\alpha)\,dt$

$$= \sum_{l=0}^{\theta(j+n)-1}\varphi_l \int_{-\infty}^{\infty} t^s[G(t)/1-\overline{\alpha}\,\overline{G}(t)]^l\,[\alpha\, g(t)/[1-\overline{\alpha}\,\overline{G}(t)]^2]\,dt$$

$$= \sum_{l=0}^{\theta(j+n)-1}\varphi_l\;\Gamma_{s,\,l,\,0}$$

Where $\Gamma_{p,q,r} = \int_{-\infty}^{\infty} t^p\,[G(t)/1-\overline{\alpha}\,\overline{G}(t)]^q\,[\alpha\,\overline{G}(t)/1-\overline{\alpha}\,\overline{G}(t)]^r\,[\alpha\, g(t)/[1-\overline{\alpha}\,\overline{G}(t)]^2]dt$

is the PWM of $\mathrm{M}O(\alpha)$ distribution.

Therefore the moments of the $\mathrm{BGMO\text{-}G}(t;\alpha,\theta,m,n)$ can be expresses in terms of the PWMs of $\mathrm{M}O(\alpha)$ (Marshall and Olkin, 1997). The PWM method can generally be used for estimating parameters quantiles of generalized distributions. These moments have low variance and no severe biases, and they compare favourably with estimators obtained by maximum likelihood.

Proceeding as above we can derive $s^{th}$ moment of the $r^{th}$ order statistic $T_{r:n}$, in a random sample of size *n* from $BGMO-G$ on using equation (17) as

$$E(T^{s}{}_{r;n}) = \sum_{l=0}^{\theta(j+n)-1} \sum_{k=0}^{\infty} \xi_{l,k}\, \Gamma_{s,k+l,\,0}$$

Where $\delta_j$, $\varphi_l$ and $\xi_{l,k}$ defined in above

**4.5 Moment generating function**

The moment generating function of $BGMO-G$ family can be easily expressed in terms of those of the exponentiated *MO* (Marshall and Olkin, 1997) distribution using the results of section 4.1. For example using equation (11) it can be seen that

$$M_T(s) = E[e^{sT}] = \int_{-\infty}^{\infty} e^{st} f^{BGMO}(t;\alpha,\theta,m,n)\,dt = \int_{-\infty}^{\infty} e^{st} \sum_{j=0}^{m-1} \delta'_j \frac{d}{dt}[\overline{F}^{MO}(t;\alpha)]^{\theta(j+n)}\, dt$$

$$= \sum_{j=0}^{m-1} \delta'_j \int_{-\infty}^{\infty} e^{st} \frac{d}{dt}[\overline{F}^{MO}(t;\alpha)]^{\theta(j+n)}\, dt = \sum_{j=0}^{m-1} \delta'_j\, M_X(s)$$

Where $M_X(s)$ is the mgf of a *MO* (Marshall and Olkin, 1997) distribution.

**4.6 Renyi Entropy**

The entropy of a random variable is a measure of uncertainty variation and has been used in various situations in science and engineering. The Rényi entropy is defined by

$$I_R(\delta) = (1-\delta)^{-1} \log\left( \int_{-\infty}^{\infty} f(t)^{\delta}\, dt \right)$$

where $\delta > 0$ and $\delta \neq 1$ For furthers details, see Song (2001). Using binomial expansion in (6) we can write

$$f^{BGMOG}(t;\alpha,\theta,m,n)^{\delta}$$

$$= \left[ \frac{\theta}{B(m,n)} \frac{\alpha\, g(t)}{[1-\overline{\alpha}\,\overline{G}(t)]^2} \left( \frac{\alpha\,\overline{G}(t)}{[1-\overline{\alpha}\,\overline{G}(t)]} \right)^{\theta-1} \left[ 1 - \left[ \frac{\alpha\,\overline{G}(t)}{1-\overline{\alpha}\,\overline{G}(t)} \right]^{\theta} \right]^{m-1} \left[ \left[ \frac{\alpha\,\overline{G}(t)}{1-\overline{\alpha}\,\overline{G}(t)} \right]^{\theta} \right]^{n-1} \right]^{\delta}$$

$$= \frac{\theta^{\delta}}{B(m,n)^{\delta}} f^{MO}(t,\alpha)^{\delta} [\overline{F}^{MO}(t,\alpha)]^{\delta(\theta-1)} [1-\{\overline{F}^{MO}(t,\alpha)\}^{\theta}]^{\delta(m-1)} [\overline{F}^{MO}(t,\alpha)]^{\theta\delta(n-1)}$$

$$= \frac{\theta^{\delta}}{B(m,n)^{\delta}} f^{MO}(t,\alpha)^{\delta} [\overline{F}^{MO}(t,\alpha)]^{\delta(\theta n-1)} [1-\{\overline{F}^{MO}(t,\alpha)\}^{\theta}]^{\delta(m-1)}$$

$$= \frac{\theta^{\delta}}{B(m,n)^{\delta}} f^{MO}(t,\alpha)^{\delta} [\overline{F}^{MO}(t,\alpha)]^{\delta(\theta n-1)} \sum_{j=0}^{\delta(m-1)} \binom{\delta(m-1)}{j} (-1)^j [\overline{F}^{MO}(t;\alpha)]^{\theta j}$$

$$= \frac{\theta^{\delta}}{B(m,n)^{\delta}} f^{MO}(t,\alpha)^{\delta} \sum_{j=0}^{\delta(m-1)} \binom{\delta(m-1)}{j} (-1)^j [\overline{F}^{MO}(t;\alpha)]^{\theta j+\delta(\theta n-1)}$$

Thus the Rényi entropy of $T$ can be obtained as

$$I_R(\delta)=(1-\delta)^{-1}\log\left(\sum_{j=0}^{\delta(m-1)} Z_j \int_{-\infty}^{\infty} f^{MO}(t,\alpha)^{\delta}\ [\overline{F}^{MO}(t;\alpha)]^{\theta j+\delta(\theta n-1)}\,dt\right)$$

Where $Z_j=\dfrac{\theta^{\delta}}{B(m,n)^{\delta}}\dbinom{\delta(m-1)}{j}(-1)^{j}$

**4.7 Quantile power series and random sample generation**

The quantile function *T*, let $t=Q(u)=F^{-1}(u)$, can be obtained by inverting (7). Let $z=Q_{m,n}(u)$ be the beta quantile function. Then,

$$t=Q(u)=Q_G\left[\frac{\alpha[1-\{1-Q_{m,n}(u)\}^{1/\theta}]}{1-\overline{\alpha}[1-\{1-Q_{m,n}(u)\}^{1/\theta}]}\right]$$

It is possible to obtain an expansion for $Q_{m,n}(u)$ in the Wolfram website as

It is possible to obtain an expansion for $Q_{m,n}(u)$ as

$$z=Q_{m,n}(u)=\sum_{i=0}^{\infty} e_i\, u^{\,i/m}$$

(see "Power series" From MathWorld--A Wolfram Web Resource. http://mathworld.wolfram.com/PowerSeries.html)

Where $e_i=[m\,B(m,n)]^{1/m}d_i$ and $d_0=0,\ d_1=1,\ d_2=(n-1)/(m+1)$,

$$d_3=\frac{(n-1)(m^2+3mn-m+5n-4)}{2(m+1)^2(m+2)}$$

$$d_4=(n-1)[m^4+(6n-1)m^3+(n+2)(8n-5)m^2+(33n^2-30n+4)m$$
$$+n(31n-47)+18]/[3(m+1)^3(m+2)(m+3)]\ldots$$

The Bowley skewness (Kenney and Keeping 1962) measures and Moors kurtosis (Moors 1988) measure are robust and less sensitive to outliers and exist even for distributions without moments. For $BGMO-G$ family these measures are given by

$$B=\frac{Q(3/4)+Q(1/4)-2Q(1/2)}{Q(3/4)-Q(1/4)}\quad\text{and}\quad M=\frac{Q(3/8)-Q(1/8)+Q(7/8)-Q(5/8)}{Q(6/8)-Q(2/8)}$$

For example, let the $G$ be exponential distribution with parameter $\lambda>0,$ having pdf and cdf as $g(t:\lambda)=\lambda e^{-\lambda t},\ t>0$ and $G(t:\lambda)=1-e^{-\lambda t}$, respectively. Then the $p^{th}$ quantile is obtained as $-(1/\lambda)\log[1-p]$. Therefore, the $p^{th}$ quantile $t_p$, of $BGMO-E$ is given by

$$t_p = -\frac{1}{\lambda}\log\left[1-\frac{\alpha[1-\{1-Q_{m,n}(p)\}^{1/\theta}]}{1-\bar{\alpha}[1-\{1-Q_{m,n}(p)\}^{1/\theta}]}\right]$$

### 4.8 Asymptotes

Here we investigate the asymptotic shapes of the proposed family following the method followed in Alizadeh *et al.*, (2015).

**Proposition 1.** The asymptotes of equations (6), (7) and (8) as $t \to 0$ are given by

$$f(t) \sim \frac{\theta\, g(t)}{B(m,n)\,\alpha}[1-\{\alpha/1-\bar{\alpha}\}^{\theta}]^{m-1} \qquad \text{as } G(t)\to 0$$

$$F(t) \sim \frac{1}{B(m,n)\,m}[1-\{\alpha/1-\bar{\alpha}\}^{\theta}]^{m} \qquad \text{as } G(t)\to 0$$

$$h(t) \sim \frac{\theta\, g(t)}{B(m,n)\,\alpha}[1-\{\alpha/1-\bar{\alpha}\}^{\theta}]^{m-1} \qquad \text{as } G(t)\to 0$$

**Proposition 2.** The asymptotes of equations (6), (7) and (8) as $t \to \infty$ are given by

$$f(t) \sim \frac{\theta\alpha^{\theta n}\, g(t)\,\bar{G}(t)^{\theta n-1}}{B(m,n)} \qquad \text{as } t\to\infty$$

$$1-F(t) \sim \frac{[\alpha\,\bar{G}(t)]^{\theta n}}{n\,B(m,n)} \qquad \text{as } t\to\infty$$

$$h(t) \sim \theta\, g(t)\bar{G}(t)^{-1} \qquad \text{as } t\to\infty$$

## 5. Estimation

### 5.1 Maximum likelihood method

The model parameters of the $BGMO-G$ distribution can be estimated by maximum likelihood. Let $t=(t_1,t_2,...t_n)^T$ be a random sample of size $n$ from $BGMO-G$ with parameter vector $\boldsymbol{\theta}=(m,n,\theta,\alpha,\boldsymbol{\beta}^T)^T$, where $\boldsymbol{\beta}=(\beta_1,\beta_2,...\beta_q)^T$ corresponds to the parameter vector of the baseline distribution $G$. Then the log-likelihood function for $\boldsymbol{\theta}$ is given by

$$\ell=\ell(\boldsymbol{\theta}) = r\log\theta + r\theta\log\alpha + \sum_{i=0}^{r}\log[g(t_i,\boldsymbol{\beta})] + (\theta-1)\sum_{i=0}^{r}\log[\bar{G}(t_i,\boldsymbol{\beta})] - r\log[B(m,n)]$$

$$-(\theta+1)\sum_{i=0}^{r}\log[1-\bar{\alpha}\,\bar{G}(t_i,\boldsymbol{\beta})] + (m-1)\sum_{i=0}^{r}\log[\,1-[\alpha\,\bar{G}(t_i,\boldsymbol{\beta})/\{1-\bar{\alpha}\,\bar{G}(t_i,\boldsymbol{\beta})\}]^{\theta}\,]$$

$$+(n-1)\sum_{i=0}^{r}\log[[\alpha\,\bar{G}(t_i,\boldsymbol{\beta})/\{1-\bar{\alpha}\,\bar{G}(t_i,\boldsymbol{\beta})\}]^{\theta}\,] \qquad (18)$$

This log-likelihood function can not be solved analytically because of its complex form but it can be maximized numerically by employing global optimization methods available with software's like R, SAS, Mathematica or by solving the nonlinear likelihood equations obtained by differentiating (18).

By taking the partial derivatives of the log-likelihood function with respect to $m, n, \theta, \alpha$ and $\boldsymbol{\beta}$ components of the score vector

$U_{\boldsymbol{\theta}} = (U_m, U_n, U_\theta, U_\alpha, U_{\beta^T})^T$ can be obtained as follows:

$$U_m = \frac{\partial \ell}{\partial m} = -r\psi(m) + r\psi(m+n) + \sum_{i=0}^{r} \log[\,1 - [\alpha \overline{G}(t_i,\boldsymbol{\beta}) / \{1 - \overline{\alpha}\,\overline{G}(t_i,\boldsymbol{\beta})\}]^\theta\,]$$

$$U_n = \frac{\partial \ell}{\partial n} = -r\psi(n) + r\psi(m+n) + b\sum_{i=0}^{r} \log[[\alpha \overline{G}(t_i,\boldsymbol{\beta}) / \{1 - \overline{\alpha}\,\overline{G}(t_i,\boldsymbol{\beta})\}]^\theta\,]$$

$$U_\theta = \frac{\partial \ell}{\partial \theta} = \frac{r}{\theta} + r\log\alpha + \sum_{i=0}^{r} \log[\overline{G}(t_i,\boldsymbol{\beta})] - \sum_{i=0}^{r} \log[1 - \overline{\alpha}\,\overline{G}(t_i,\boldsymbol{\beta})]$$

$$+(1-m)\sum_{i=0}^{r} \frac{[\alpha \overline{G}(t_i,\boldsymbol{\beta}) / \{1 - \overline{\alpha}\,\overline{G}(t_i,\boldsymbol{\beta})\}]^\theta \log[\alpha \overline{G}(t_i,\boldsymbol{\beta}) / \{1 - \overline{\alpha}\,\overline{G}(t_i,\boldsymbol{\beta})\}]}{1 - [\alpha \overline{G}(t_i,\boldsymbol{\beta}) / \{1 - \overline{\alpha}\,\overline{G}(t_i,\boldsymbol{\beta})\}]^\theta}$$

$$+(n-1)\sum_{i=0}^{r} \frac{[\alpha \overline{G}(t_i,\boldsymbol{\beta}) / \{1 - \overline{\alpha}\,\overline{G}(t_i,\boldsymbol{\beta})\}]^\theta \log[\alpha \overline{G}(t_i,\boldsymbol{\beta}) / \{1 - \overline{\alpha}\,\overline{G}(t_i,\boldsymbol{\beta})\}]}{[\alpha \overline{G}(t_i,\boldsymbol{\beta}) / \{1 - \overline{\alpha}\,\overline{G}(t_i,\boldsymbol{\beta})\}]^\theta}$$

$$U_\alpha = \frac{\partial \ell}{\partial \alpha} = \frac{r\theta}{\alpha} - (\theta+1)\sum_{i=0}^{r} \frac{\overline{G}(t_i,\boldsymbol{\beta})}{1 - \overline{\alpha}\,\overline{G}(t_i,\boldsymbol{\beta})}$$

$$+\theta\alpha(1-m)\sum_{i=0}^{r} \frac{\overline{G}(t_i,\boldsymbol{\beta})^\theta G(t_i,\boldsymbol{\beta})}{[\{1 - \overline{\alpha}\,\overline{G}(t_i,\boldsymbol{\beta})\}^\theta - \{\alpha \overline{G}(t_i,\boldsymbol{\beta})\}^\theta][1 - \overline{\alpha}\,\overline{G}(t_i,\boldsymbol{\beta})]} + \frac{1}{\alpha}\theta(n-1)\sum_{i=0}^{r} \frac{G(t_i,\boldsymbol{\beta})}{1 - \overline{\alpha}\,\overline{G}(t_i,\boldsymbol{\beta})}$$

$$U_\beta = \frac{\partial \ell}{\partial \boldsymbol{\beta}} = \sum_{i=0}^{r} \frac{g^{(\boldsymbol{\beta})}(t_i,\boldsymbol{\beta})}{g(t_i,\boldsymbol{\beta})} + (1-\theta)\sum_{i=0}^{r} \frac{G^{(\boldsymbol{\beta})}(t_i,\boldsymbol{\beta})}{\overline{G}(t_i,\boldsymbol{\beta})} - (\theta+1)\sum_{i=0}^{r} \frac{\overline{\alpha}\,G^{(\boldsymbol{\beta})}(t_i,\boldsymbol{\beta})}{1 - \overline{\alpha}\,\overline{G}(t_i,\boldsymbol{\beta})}$$

$$+\theta\alpha^\theta(m-1)\sum_{i=0}^{r} \frac{\overline{G}(t_i,\boldsymbol{\beta})^{\theta-1} G^{(\boldsymbol{\beta})}(t_i,\boldsymbol{\beta})}{[\{1 - \overline{\alpha}\,\overline{G}(t_i,\boldsymbol{\beta})\}^\theta - \{\alpha \overline{G}(t_i,\boldsymbol{\beta})\}^\theta][1 - \overline{\alpha}\,\overline{G}(t_i,\boldsymbol{\beta})]}$$

$$+\theta(1-n)\sum_{i=0}^{r} \frac{G^{(\boldsymbol{\beta})}(t_i,\boldsymbol{\beta})}{[1 - \overline{\alpha}\,\overline{G}(t_i,\boldsymbol{\beta})]\overline{G}(t_i,\boldsymbol{\beta})}$$

Where $\psi(.)$ is the digamma function.

**5.2 Asymptotic standard error and confidence interval for the mles:**

The asymptotic variance-covariance matrix of the MLEs of parameters can obtained by inverting the Fisher information matrix $\mathrm{I}(\boldsymbol{\theta})$ which can be derived using the second partial derivatives of the log-likelihood function with respect to each parameter. The $i\,j^{th}$ elements of $\mathrm{I}_n(\boldsymbol{\theta})$ are given by

$$\mathrm{I}_{ij} = -E\left(\frac{\partial^2 l(\boldsymbol{\theta})}{\partial\theta_i\,\partial\theta_j}\right), \qquad i,\ j = 1, 2, \cdots, 3+q$$

The exact evaluation of the above expectations may be cumbersome. In practice one can estimate $\mathrm{I}_n(\boldsymbol{\theta})$ by the observed Fisher's information matrix $\hat{\mathrm{I}}_n(\hat{\boldsymbol{\theta}})$ is defined as:

$$\hat{\mathrm{I}}_{ij} \approx \left( -\frac{\partial^2 l(\boldsymbol{\theta})}{\partial \theta_i \partial \theta_j} \right)_{\boldsymbol{\theta}=\hat{\boldsymbol{\theta}}}, \qquad i,\ j = 1, 2, \cdots, 3+q$$

Using the general theory of MLEs under some regularity conditions on the parameters as $n \to \infty$ the asymptotic distribution of $\sqrt{n}(\hat{\boldsymbol{\theta}} - \boldsymbol{\theta})$ is $N_k(0, V_n)$ where $V_n = (v_{jj}) = \mathrm{I}_n^{-1}(\boldsymbol{\theta})$. The asymptotic behaviour remains valid if $V_n$ is replaced by $\hat{V}_n = \hat{\mathrm{I}}^{-1}(\hat{\boldsymbol{\theta}})$. This result can be used to provide large sample standard errors and also construct confidence intervals for the model parameters. Thus an approximate standard error and $(1-\gamma/2)100\%$ confidence interval for the mle of $j^{\text{th}}$ parameter $\theta_j$ are respectively given by $\sqrt{\hat{v}_{jj}}$ and $\hat{\theta}_j \pm Z_{\gamma/2}\sqrt{\hat{v}_{jj}}$ , where $Z_{\gamma/2}$ is the $\gamma/2$ point of standard normal distribution.

As an illustration on the MLE method its large sample standard errors, confidence interval in the case of $BGMO - E(m, n, \theta, \alpha, \lambda)$ is discussed in an appendix.

**5.3 Real life applications**

In this subsection, we consider fitting of three real data sets to show that the proposed $BGMO - G$ distribution can be a better model than $GMOKw - G$ (Handique and Chakraborty, 2015) by taking as Weibull distribution as *G*. We have estimated the parameters by numerical maximization of loglikelihood function and provided their standard errors and 95% confidence intervals using large sample approach (see appendix).

In order to compare the distributions, we have considered known criteria like AIC (Akaike Information Criterion), BIC (Bayesian Information Criterion), CAIC (Consistent Akaike Information Criterion) and HQIC (Hannan-Quinn Information Criterion). It may be noted that $AIC = 2k - 2l$;

$$BIC = k\log(n) - 2l;\ CAIC = AIC + \frac{2k(k+1)}{n-k-1};\ \text{and}\ \ HQIC = 2k\log[\log(n)] - 2l$$

where $k$ is the number of parameters in the statistical model, $n$ the sample size and $l$ is the maximized value of the log-likelihood function under the considered model. In these applications method of maximum likelihood will be used to obtain the estimate of parameters.

**Example I:**

The following data set gives the time to failure $(10^3 h)$ of turbocharger of one type of engine given in Xu et al. (2003).

{1.6, 2.0, 2.6, 3.0, 3.5, 3.9, 4.5, 4.6, 4.8, 5.0, 5.1, 5.3, 5.4, 5.6, 5.8, 6.0, 6.0, 6.1, 6.3, 6.5, 6.5, 6.7, 7.0, 7.1, 7.3, 7.3, 7.3, 7.7, 7.7, 7.8, 7.9, 8.0, 8.1, 8.3, 8.4, 8.4, 8.5, 8.7, 8.8, 9.0}

**Table 1:** MLEs, standard errors and 95% confidence intervals (in parentheses) and the AIC, BIC, CAIC and HQIC values for the data set.

| Parameters | $GMOKw-W$ | $BGMO-W$ |
|---|---|---|
| $\hat{a}$ | 1.178<br>(0.017)<br>(1.14, 1.21) | 1.187<br>(0.702)<br>(-0.19, 2.56) |
| $\hat{b}$ | 0.291<br>(0.209)<br>(-0.12, 0.70) | 2.057<br>(2.240)<br>(-2.33, 6.45) |
| $\hat{\lambda}$ | 0.617<br>(0.002)<br>(0.61, 0.62) | 0.009<br>(0.006)<br>(-0.00276, 0.02) |
| $\hat{\beta}$ | 1.855<br>(0.003)<br>(1.85, 1.86) | 4.194<br>(0.668)<br>(2.88, 5.50) |
| $\hat{\alpha}$ | 1.619<br>(0.977)<br>(-0.29, 3.53) | 0.047<br>(0.108)<br>(-0.16, 0.26) |
| $\hat{\theta}$ | 0.178<br>(0.144)<br>(-0.10, 0.46) | 0.017<br>(0.016)<br>(-0.01, 0.05) |
| log-likelihood( $l_{max}$ ) | -90.99 | **-80.38** |
| AIC | 193.98 | **172.76** |
| BIC | 204.11 | **182.89** |
| CAIC | 196.53 | **175.31** |
| HQIC | 197.65 | **176.43** |

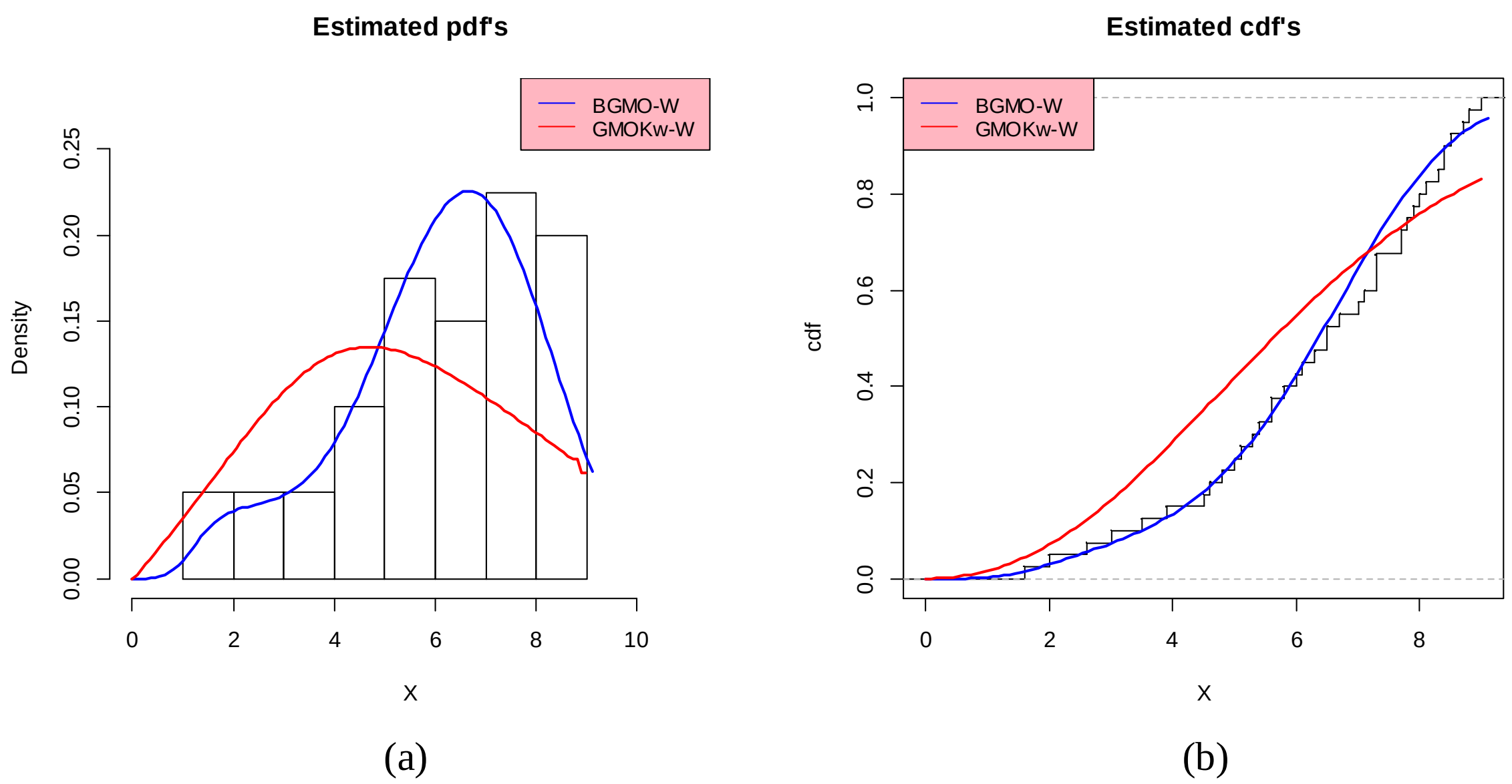

**Fig: 3** Plots of the observed (a) histogram and estimated pdf's and (b) estimated cdf's for the $BGMO-G$ and $GMOKw-G$ for example I.

**Example II:**

The following data is about 346 nicotine measurements made from several brands of cigarettes in 1998. The data have been collected by the Federal Trade Commission which is an independent agency of the US government, whose main mission is the promotion of consumer protection. **[http: //www.ftc.gov/ reports/tobacco or** http**: // pw1.netcom.com/ rdavis2/ smoke. html.]**

{1.3, 1.0, 1.2, 0.9, 1.1, 0.8, 0.5, 1.0, 0.7, 0.5, 1.7, 1.1, 0.8, 0.5, 1.2, 0.8, 1.1, 0.9, 1.2, 0.9, 0.8, 0.6, 0.3, 0.8, 0.6, 0.4, 1.1, 1.1, 0.2, 0.8, 0.5, 1.1, 0.1, 0.8, 1.7, 1.0, 0.8, 1.0, 0.8, 1.0, 0.2, 0.8, 0.4, 1.0, 0.2, 0.8, 1.4, 0.8, 0.5, 1.1, 0.9, 1.3, 0.9, 0.4, 1.4, 0.9, 0.5, 1.7, 0.9, 0.8, 0.8, 1.2, 0.9, 0.8, 0.5, 1.0, 0.6, 0.1, 0.2, 0.5, 0.1, 0.1, 0.9, 0.6, 0.9, 0.6, 1.2, 1.5, 1.1, 1.4, 1.2, 1.7, 1.4, 1.0, 0.7, 0.4, 0.9, 0.7, 0.8, 0.7, 0.4, 0.9, 0.6, 0.4, 1.2, 2.0, 0.7, 0.5, 0.9, 0.5, 0.9, 0.7, 0.9, 0.7, 0.4, 1.0, 0.7, 0.9, 0.7, 0.5, 1.3, 0.9, 0.8, 1.0, 0.7, 0.7, 0.6, 0.8, 1.1, 0.9, 0.9, 0.8, 0.8, 0.7, 0.7, 0.4, 0.5, 0.4, 0.9, 0.9, 0.7, 1.0, 1.0, 0.7, 1.3, 1.0, 1.1, 1.1, 0.9, 1.1, 0.8, 1.0, 0.7, 1.6, 0.8, 0.6, 0.8, 0.6, 1.2, 0.9, 0.6, 0.8, 1.0, 0.5, 0.8, 1.0, 1.1, 0.8, 0.8, 0.5, 1.1, 0.8, 0.9, 1.1, 0.8, 1.2, 1.1, 1.2, 1.1, 1.2, 0.2, 0.5, 0.7, 0.2, 0.5, 0.6, 0.1, 0.4, 0.6, 0.2, 0.5, 1.1, 0.8, 0.6, 1.1, 0.9, 0.6, 0.3, 0.9, 0.8, 0.8, 0.6, 0.4, 1.2, 1.3, 1.0, 0.6, 1.2, 0.9, 1.2, 0.9, 0.5, 0.8, 1.0, 0.7, 0.9, 1.0, 0.1, 0.2, 0.1, 0.1, 1.1, 1.0, 1.1, 0.7, 1.1, 0.7, 1.8, 1.2, 0.9, 1.7, 1.2, 1.3, 1.2, 0.9, 0.7, 0.7, 1.2, 1.0, 0.9, 1.6, 0.8, 0.8, 1.1, 1.1, 0.8, 0.6, 1.0, 0.8, 1.1, 0.8, 0.5, 1.5, 1.1, 0.8, 0.6, 1.1, 0.8, 1.1, 0.8, 1.5, 1.1, 0.8, 0.4, 1.0, 0.8, 1.4, 0.9, 0.9, 1.0, 0.9, 1.3, 0.8, 1.0, 0.5, 1.0, 0.7, 0.5, 1.4, 1.2, 0.9, 1.1, 0.9, 1.1, 1.0, 0.9, 1.2, 0.9, 1.2, 0.9, 0.5, 0.9, 0.7, 0.3, 1.0, 0.6, 1.0, 0.9, 1.0, 1.1, 0.8, 0.5, 1.1, 0.8, 1.2, 0.8, 0.5, 1.5, 1.5, 1.0, 0.8, 1.0, 0.5, 1.7, 0.3, 0.6, 0.6, 0.4, 0.5, 0.5, 0.7, 0.4, 0.5, 0.8, 0.5, 1.3, 0.9, 1.3, 0.9, 0.5, 1.2, 0.9, 1.1, 0.9, 0.5, 0.7, 0.5, 1.1, 1.1, 0.5, 0.8, 0.6, 1.2, 0.8, 0.4, 1.3, 0.8, 0.5, 1.2, 0.7, 0.5, 0.9, 1.3, 0.8, 1.2, 0.9}

**Table 2:** MLEs, standard errors and 95% confidence intervals (in parentheses) and the AIC, BIC, CAIC and HQIC values for the nicotine measurements data.

| Parameters | $GMOKw-W$ | $BGMO-W$ |
|---|---|---|
| $\hat{a}$ | 0.765<br>(0.025)<br>(0.72, 0.81) | 0.866<br>(0.159)<br>(0.55, 1.18) |
| $\hat{b}$ | 2.139<br>(0.774)<br>(0.62, 3.66) | 0.329<br>(0.167)<br>(0.00168, 0.66) |
| $\hat{\lambda}$ | 4.271<br>(0.018)<br>(4.24, 4.31) | 2.131<br>(1.182)<br>(-0.19, 4.45) |
| $\hat{\beta}$ | 2.919<br>(0.013)<br>(2.89, 2.94) | 2.223<br>(0.725)<br>(0.80, 3.64) |
| $\hat{\alpha}$ | 1.097<br>(0.309)<br>(0.49, 1.70) | 3.285<br>(3.901)<br>(-4.36, 10.93) |
| $\hat{\theta}$ | 0.114<br>(0.042)<br>(0.03, 0.19) | 2.635<br>(1.310)<br>(0.07, 5.20) |
| log-likelihood( $l_{\max}$ ) | -111.75 | **-109.28** |
| AIC | 235.50 | **230.56** |
| BIC | 258.58 | **253.64** |
| CAIC | 235.75 | **230.80** |
| HQIC | 244.69 | **239.76** |

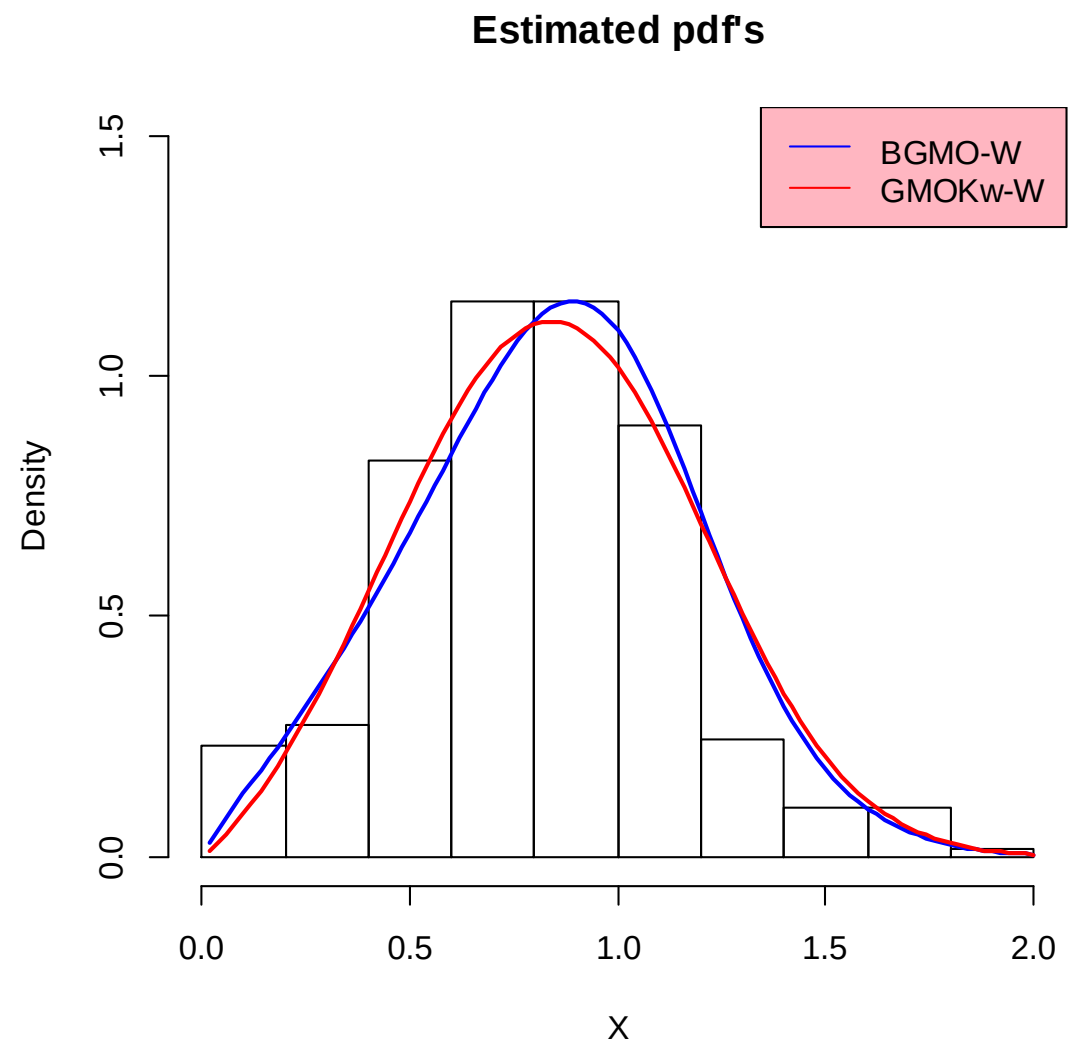

(a)

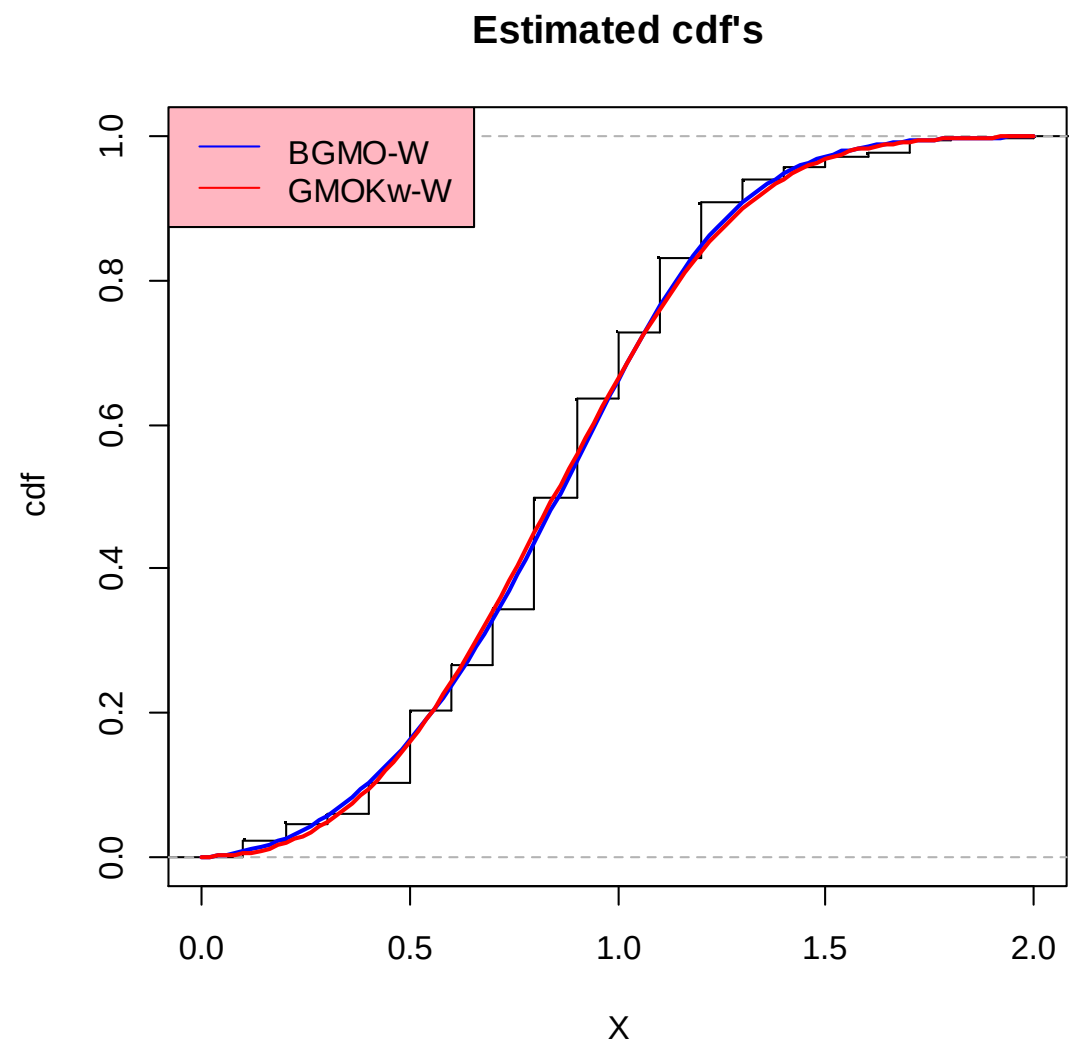

(b)

**Fig: 4** Plots of the observed (a) histogram and estimated pdf's and (b) estimated cdf's for the $BGMO-G$ and $GMOKw-G$ for example II.

**Example III:**

This data set consists of 100 observations of breaking stress of carbon fibres (in Gba) given by Nichols and Padgett (2006).

{3.70, 2.74, 2.73, 2.50, 3.60, 3.11, 3.27, 2.87, 1.47, 3.11, 4.42, 2.40, 3.15, 2.67,3.31, 2.81, 0.98, 5.56, 5.08, 0.39, 1.57, 3.19, 4.90, 2.93, 2.85, 2.77, 2.76, 1.73, 2.48, 3.68, 1.08, 3.22, 3.75, 3.22, 2.56, 2.17, 4.91, 1.59, 1.18, 2.48, 2.03, 1.69, 2.43, 3.39, 3.56, 2.83, 3.68, 2.00, 3.51, 0.85, 1.61, 3.28, 2.95, 2.81, 3.15, 1.92, 1.84, 1.22, 2.17, 1.61, 2.12, 3.09, 2.97, 4.20, 2.35, 1.41, 1.59, 1.12, 1.69, 2.79, 1.89, 1.87, 3.39, 3.33, 2.55, 3.68, 3.19, 1.71, 1.25, 4.70, 2.88, 2.96, 2.55, 2.59, 2.97, 1.57, 2.17, 4.38, 2.03, 2.82, 2.53, 3.31, 2.38, 1.36, 0.81, 1.17, 1.84, 1.80, 2.05, 3.65}.

**Table 3:** MLEs, standard errors and 95% confidence intervals (in parentheses) and the AIC, BIC, CAIC and HQIC values for the breaking stress of carbon fibres data.

| Parameters | $GMOKw-W$ | $BGMO-W$ |
|---|---|---|
| $\hat{a}$ | 1.015<br>(0.071)<br>(0.88, 1.15) | 1.458<br>(1.123)<br>(-0.74, 3.66) |
| $\hat{b}$ | 0.385<br>(0.168)<br>(0.06, 0.71) | 0.734<br>(1.810)<br>(-2.81, 4.28) |
| $\hat{\lambda}$ | 0.803<br>(0.003)<br>(0.79, 0.81) | 0.598<br>(1.496)<br>(-2.33, 3.53) |
| $\hat{\beta}$ | 2.222<br>(0.004)<br>(2.21, 2.23) | 2.439<br>(0.779)<br>(0.91, 3.97) |
| $\hat{\alpha}$ | 1.482<br>(0.440)<br>(0.62, 2.34) | 0.685<br>(1.150)<br>(-1.57, 2.94) |
| $\hat{\theta}$ | 0.345<br>(0.169)<br>(0.01, 0.68) | 0.201<br>(0.573)<br>(-0.92, 1.32) |
| log-likelihood( $l_{\max}$ ) | 142.63 | **-141.29** |
| AIC | 297.26 | **294.58** |
| BIC | 312.89 | **310.21** |
| CAIC | 298.16 | **295.48** |
| HQIC | 303.59 | **300.92** |

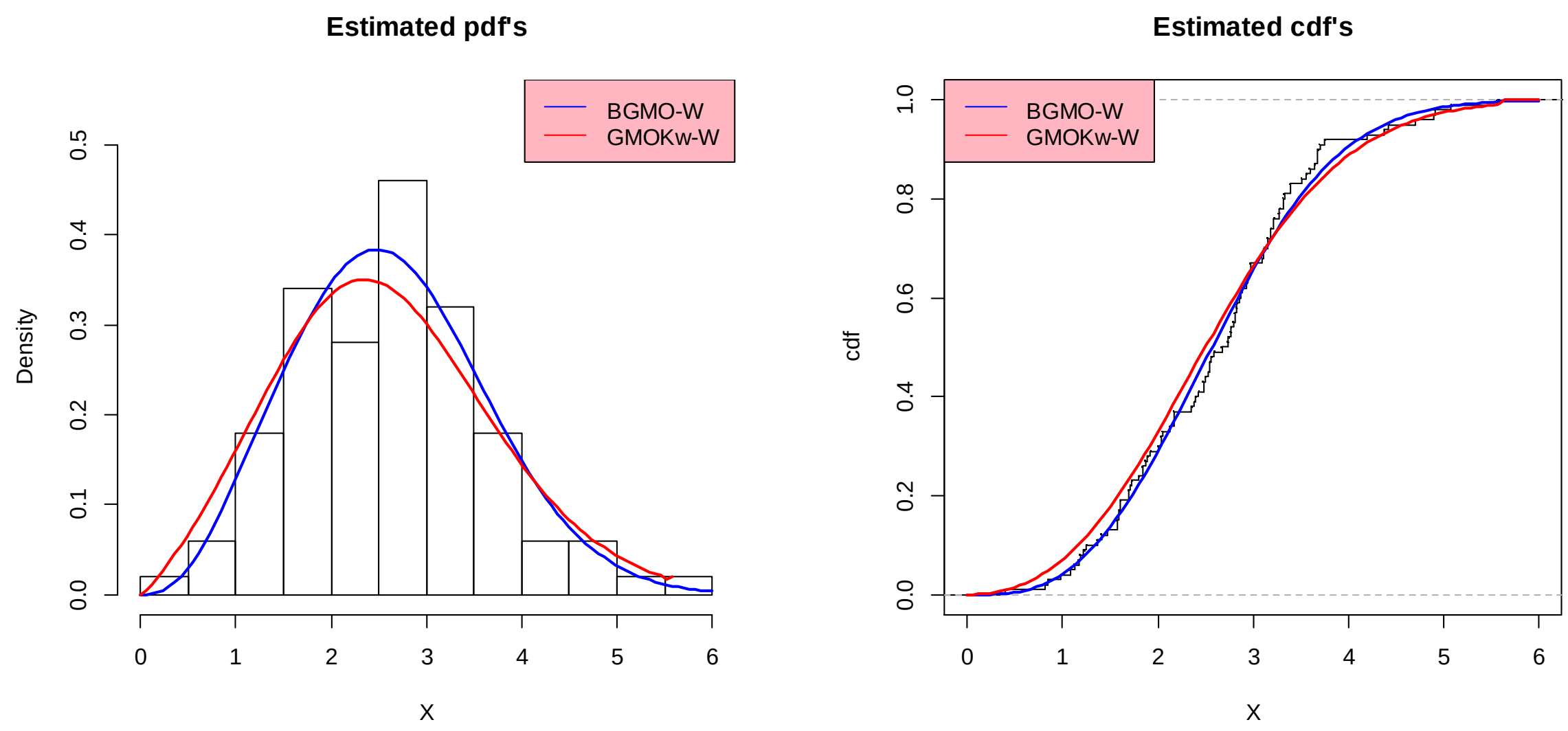

(a) (b)

**Fig: 5** Plots of the observed (a) histogram and estimated pdf's and (b) estimated cdf's for the $BGMO-G$ and $GMOKw-G$ for example III.

In Tables 1, 2 and 3 the MLEs, se's (in parentheses) and 95% confidence intervals (in parentheses) of the parameters for the fitted distributions along with AIC, BIC, CAIC and HQIC values are presented for example I, II, and III respectively. In the entire examples considered here based on the lowest values of the AIC, BIC, CAIC and HQIC, the $BGMO-W$ distribution turns out to be a better distribution than $GMOKw-W$ distribution. A visual comparison of the closeness of the fitted densities with the observed histogram of the data for example I, II, and III are presented in the figures 3, 4, and 5 respectively. These plots indicate that the proposed distributions provide a closer fit to these data.

## 6. Conclusion

Beta Generalized Marshall-Olkin family of distributions is introduced and some of its important properties are studied. The maximum likelihood and moment method for estimating the parameters are also discussed. Application of three real life data fitting shows good result in favour of the proposed family when compared to Generalized Marshall-Olkin Kumaraswamy extended family of distributions. It is therefore expected that this family of distribution will be an important addition to the existing literature on distribution.

## Appendix: Maximum likelihood estimation for $BGMO-E$

The pdf of the $BGMO-E$ distribution is given by

$$f^{BGMOE}(t)=\frac{1}{B(m,n)}\frac{\theta\alpha^{\theta}\lambda e^{-\lambda t}[e^{-\lambda t}]^{\theta-1}}{[1-\bar{\alpha}e^{-\lambda t}]^{\theta+1}}\left[1-\left[\frac{\alpha e^{-\lambda t}}{1-\bar{\alpha}e^{-\lambda t}}\right]^{\theta}\right]^{m-1}\left[\left[\frac{\alpha e^{-\lambda t}}{1-\bar{\alpha}e^{-\lambda t}}\right]^{\theta}\right]^{n-1}$$

For a random sample of size $n$ from this distribution, the log-likelihood function for the parameter vector $\boldsymbol{\theta} = (m,n,\theta,\alpha,\lambda)^T$ is given by

$$\ell = \ell(\boldsymbol{\theta}) = r\log\theta + r\theta\log\alpha + r\log\lambda - \lambda\sum_{i=0}^{r} t_i + (\theta-1)\sum_{i=0}^{r}\log(e^{-\lambda t_i}) - r\log[B(m,n)]$$

$$-(\theta+1)\sum_{i=0}^{r}\log[1-\bar{\alpha}\,e^{-\lambda t_i}] + (m-1)\sum_{i=0}^{r}\log[1-[\,\alpha\,e^{-\lambda t_i}\big/\{1-\bar{\alpha}\,e^{-\lambda t_i}\}]^{\theta}\,]$$

$$+(n-1)\sum_{i=0}^{r}\log\,[\,\alpha\,e^{-\lambda t_i}\big/\{1-\bar{\alpha}\,e^{-\lambda t_i}\}]^{\theta}$$

The components of the score vector $\boldsymbol{\theta} = (m,n,\theta,\alpha,\lambda)^T$ are

$$U_m = \frac{\partial\ell}{\partial m} = -r\psi(m) + r\psi(m+n) + \sum_{i=0}^{r}\log[1-[\,\alpha\,e^{-\lambda t_i}\big/\{1-\bar{\alpha}\,e^{-\lambda t_i}\}]^{\theta}\,]$$

$$U_n = \frac{\partial\ell}{\partial n} = -r\psi(n) + r\psi(m+n) + \sum_{i=0}^{r}\log\,[\,\alpha\,e^{-\lambda t_i}\big/\{1-\bar{\alpha}\,e^{-\lambda t_i}\}]^{\theta}$$

$$U_\theta = \frac{\partial\ell}{\partial\theta} = \frac{r}{\theta} + r\log\alpha + \sum_{i=0}^{r}\log[e^{-\lambda t_i}] - \sum_{i=0}^{r}\log[1-\bar{\alpha}\,e^{-\lambda t_i}]$$

$$+(1-m)\sum_{i=0}^{r}\frac{[\alpha\,e^{-\lambda t_i}\big/\{1-\bar{\alpha}\,e^{-\lambda t_i}\}]^{\theta}\log[\alpha\,e^{-\lambda t_i}\big/\{1-\bar{\alpha}\,e^{-\lambda t_i}\}]}{1-[\alpha\,e^{-\lambda t_i}\big/\{1-\bar{\alpha}\,e^{-\lambda t_i}\}]^{\theta}}$$

$$+(n-1)\sum_{i=0}^{r}\frac{[\,\alpha\,e^{-\lambda t_i}\big/\{1-\bar{\alpha}\,e^{-\lambda t_i}\}]^{\theta}\log[\,\alpha\,e^{-\lambda t_i}\big/\{1-\bar{\alpha}\,e^{-\lambda t_i}\}]}{[\,\alpha\,e^{-\lambda t_i}\big/\{1-\bar{\alpha}\,e^{-\lambda t_i}\}]^{\theta}}$$

$$U_\alpha = \frac{\partial\ell}{\partial\alpha} = \frac{r\theta}{\alpha} - (\theta+1)\sum_{i=0}^{r}\frac{e^{-\lambda t_i}}{1-\bar{\alpha}\,e^{-\lambda t_i}}$$

$$+\theta\alpha(1-m)\sum_{i=0}^{r}\frac{[e^{-\lambda t_i}]^{\theta}\,[1-e^{-\lambda t_i}]}{[\{1-\bar{\alpha}\,e^{-\lambda t_i}\}^{\theta} - \{\alpha\,e^{-\lambda t_i}\}^{\theta}\,][1-\bar{\alpha}\,e^{-\lambda t_i}]} + \frac{1}{\alpha}\theta(n-1)\sum_{i=0}^{r}\frac{1-e^{-\lambda t_i}}{1-\bar{\alpha}\,e^{-\lambda t_i}}$$

$$U_\lambda = \frac{\partial\ell}{\partial\lambda} = \frac{r}{\lambda} - \sum_{i=0}^{r} t_i + (\theta-1)\sum_{i=0}^{r}\frac{\lambda\,e^{-\lambda t_i}}{1-e^{-\lambda t_i}} - (\theta+1)\sum_{i=0}^{r}\frac{\bar{\alpha}\,\lambda\,e^{-\lambda t_i}}{1-\bar{\alpha}\,e^{-\lambda t_i}}$$

$$+\theta\alpha^{\theta}\lambda(m-1)\sum_{i=0}^{r}\frac{[e^{-\lambda t_i}]^{\theta}}{[\{1-\bar{\alpha}\,e^{-\lambda t_i}\}^{\theta} - \{\alpha\,e^{-\lambda t_i}\}^{\theta}\,][1-\bar{\alpha}\,e^{-\lambda t_i}]}$$

$$+\theta\lambda(1-n)\sum_{i=0}^{r}\frac{e^{-\lambda t_i}}{[1-\bar{\alpha}\,e^{-\lambda t_i}]\,e^{-\lambda t_i}}$$

The asymptotic variance covariance matrix for mles of the unknown parameters $\boldsymbol{\theta} = (m,n,\theta,\alpha,\lambda)$ of $BGMO-E\ \ (m,n,\theta,\alpha,\lambda)$ distribution is estimated by

$$\hat{\mathrm{I}}_n{}^{-1}(\hat{\boldsymbol{\theta}}) = \begin{pmatrix} \mathrm{var}(\hat{m}) & \mathrm{cov}(\hat{m},\hat{n}) & \mathrm{cov}(\hat{m},\hat{\theta}) & \mathrm{cov}(\hat{m},\hat{\alpha}) & \mathrm{cov}(\hat{m},\hat{\lambda}) \\ \mathrm{cov}(\hat{n},\hat{m}) & \mathrm{var}(\hat{n}) & \mathrm{cov}(\hat{n},\hat{\theta}) & \mathrm{cov}(\hat{n},\hat{\alpha}) & \mathrm{cov}(\hat{n},\hat{\lambda}) \\ \mathrm{cov}(\hat{\theta},\hat{m}) & \mathrm{cov}(\hat{\theta},\hat{n}) & \mathrm{var}(\hat{\theta}) & \mathrm{cov}(\hat{\theta},\hat{\alpha}) & \mathrm{cov}(\hat{\theta},\hat{\lambda}) \\ \mathrm{cov}(\hat{\alpha},\hat{m}) & \mathrm{cov}(\hat{\alpha},\hat{n}) & \mathrm{cov}(\hat{\alpha},\hat{\theta}) & \mathrm{var}(\hat{\alpha}) & \mathrm{cov}(\hat{\alpha},\hat{\lambda}) \\ \mathrm{cov}(\hat{\lambda},\hat{m}) & \mathrm{cov}(\hat{\lambda},\hat{n}) & \mathrm{cov}(\hat{\lambda},\hat{\theta}) & \mathrm{cov}(\hat{\lambda},\hat{\alpha}) & \mathrm{var}(\hat{\lambda}) \end{pmatrix}$$

Where the elements of the information matrix $\hat{\mathbf{I}}_n(\hat{\boldsymbol{\theta}}) = \left( -\left( \frac{\partial^2 l(\boldsymbol{\theta})}{\partial\theta_i \partial\theta_j} \right) \right)_{\boldsymbol{\theta}=\hat{\boldsymbol{\theta}}}$ can be derived using the following second partial derivatives:

$$\frac{\partial^2\ell}{\partial m^2} = -r\psi'(m) - r\psi'(m+n)$$

$$\frac{\partial^2\ell}{\partial n^2} = -r\psi'(n) - r\psi'(m+n)$$

$$\frac{\partial^2\ell}{\partial\theta^2} = -\frac{r}{\theta^2} + (1-m)\sum_{i=0}^{r} \frac{[\alpha e^{-\lambda t_i} / \{1-\overline{\alpha} e^{-\lambda t_i}\}]^{2\theta} \log[\alpha e^{-\lambda t_i} / \{1-\overline{\alpha} e^{-\lambda t_i}\}]^2}{[1-[\alpha e^{-\lambda t_i} / \{1-\overline{\alpha} e^{-\lambda t_i}\}]^{\theta}]^2}$$

$$+ (1-m)\sum_{i=0}^{r} \frac{[\alpha e^{-\lambda t_i} / \{1-\overline{\alpha} e^{-\lambda t_i}\}]^{\theta} \log[\alpha e^{-\lambda t_i} / \{1-\overline{\alpha} e^{-\lambda t_i}\}]^2}{1-[\alpha e^{-\lambda t_i} / \{1-\overline{\alpha} e^{-\lambda t_i}\}]^{\theta}}$$

$$\frac{\partial^2\ell}{\partial\alpha^2} = \frac{r\theta}{\alpha^2} - (\theta+1)\sum_{i=0}^{r} \frac{e^{-2\lambda t_i}}{(1-\overline{\alpha} e^{-\lambda t_i})^2}$$

$$+ \theta\,(n-1)\sum_{i=0}^{r} \frac{e^{\lambda t_i}\,[1-\overline{\alpha} e^{-\lambda t_i}]\left(-2\alpha\, e^{-3\lambda t_i} t_i / \{1-\overline{\alpha} e^{-\lambda t_i}\}^3 - 2\alpha\, e^{-2\lambda t_i} t_i / \{1-\overline{\alpha} e^{-\lambda t_i}\}^2\right)}{b}$$

$$+ \theta\,(n-1)\sum_{i=0}^{r} \frac{-\alpha\, e^{-2\lambda t_i} t_i / \{1-\overline{\alpha} e^{-\lambda t_i}\}^2 + \alpha\, e^{-\lambda t_i} t_i / \{1-\overline{\alpha} e^{-\lambda t_i}\}}{b}$$

$$- \theta\,(n-1)\sum_{i=0}^{r} \frac{e^{\lambda t_i}\,[1-\overline{\alpha} e^{-\lambda t_i}]\left(-\alpha\, e^{-2\lambda t_i} t_i / \{1-\overline{\alpha} e^{-\lambda t_i}\}^2 + \alpha\, e^{-\lambda t_i} t_i / \{1-\overline{\alpha} e^{-\lambda t_i}\}\right)}{b^2}$$

$$+ (1-m)\sum_{i=0}^{r} \frac{\theta[\alpha e^{-\lambda t_i} / \{1-\overline{\alpha} e^{-\lambda t_i}\}]^{\theta-1}\left(2\alpha\, e^{-3\lambda t_i} / \{1-\overline{\alpha} e^{-\lambda t_i}\}^3 - 2\, e^{-2\lambda t_i} / \{1-\overline{\alpha} e^{-\lambda t_i}\}^2\right)}{1-[\alpha e^{-\lambda t_i} / \{1-\overline{\alpha} e^{-\lambda t_i}\}]^{\theta}}$$

$$+ (1-m)\sum_{i=0}^{r} \frac{\theta^2[\alpha e^{-\lambda t_i} / \{1-\overline{\alpha} e^{-\lambda t_i}\}]^{2\theta-2}\left(-\alpha\, e^{-2\lambda t_i} / \{1-\overline{\alpha} e^{-\lambda t_i}\}^2 + e^{-\lambda t_i} / \{1-\overline{\alpha} e^{-\lambda t_i}\}\right)^2}{[1-[\alpha e^{-\lambda t_i} / \{1-\overline{\alpha} e^{-\lambda t_i}\}]^{\theta}]^2}$$

$$+ (1-m)\sum_{i=0}^{r} \frac{\theta(\theta-1)[\alpha e^{-\lambda t_i} / \{1-\overline{\alpha} e^{-\lambda t_i}\}]^{\theta-2}\left(-\alpha\, e^{-2\lambda t_i} / \{1-\overline{\alpha} e^{-\lambda t_i}\}^2 + e^{-\lambda t_i} / \{1-\overline{\alpha} e^{-\lambda t_i}\}\right)^2}{1-[\alpha e^{-\lambda t_i} / \{1-\overline{\alpha} e^{-\lambda t_i}\}]^{\theta}}$$

$$\frac{\partial^2 \ell}{\partial \lambda^2} = -\frac{r}{\lambda^2} + (\theta+1)\sum_{i=0}^{r}\left(-\frac{\bar{\alpha}^2 \lambda e^{-2\lambda t_i} t_i^2}{[1-\bar{\alpha}e^{-\lambda t_i}]^2} - \frac{\bar{\alpha}\lambda e^{-\lambda t_i} t_i^2}{1-\bar{\alpha}e^{-\lambda t_i}}\right)$$

$$+\theta(n-1)\sum_{i=0}^{r}\frac{\bar{\alpha}\, t_i \left(-\bar{\alpha}\alpha e^{-2\lambda t_i} t_i \big/ \{1-\bar{\alpha}e^{-\lambda t_i}\}^2 - e^{-\lambda t_i} t_i \big/ \{1-\bar{\alpha}e^{-\lambda t_i}\}\right)}{b}$$

$$+\theta(n-1)\sum_{i=0}^{r}\frac{e^{-\lambda t_i}\, [1-\bar{\alpha}e^{-\lambda t_i}] t_i \left(-\bar{\alpha}\alpha e^{-2\lambda t_i} t_i \big/ \{1-\bar{\alpha}e^{-\lambda t_i}\}^2 - e^{-\lambda t_i} t_i \big/ \{1-\bar{\alpha}e^{-\lambda t_i}\}\right)}{b}$$

$$+\theta(n-1)\sum_{i=0}^{r}\frac{e^{-\lambda t_i}\, [1-\bar{\alpha}e^{-\lambda t_i}] t_i \begin{pmatrix} -2\bar{\alpha}^2\alpha e^{-3\lambda t_i} t_i^2 \big/ \{1-\bar{\alpha}e^{-\lambda t_i}\}^3 + 3\bar{\alpha}\alpha e^{-2\lambda t_i} t_i^2 \big/ \{1-\bar{\alpha}e^{-\lambda t_i}\}^2 \\ +\alpha e^{-\lambda t_i} t_i^2 \big/ \{1-\bar{\alpha}e^{-\lambda t_i}\} \end{pmatrix}}{b}$$

$$+(1-m)\sum_{i=0}^{r}\frac{\theta^2 [\alpha e^{-\lambda t_i} \big/ \{1-\bar{\alpha}e^{-\lambda t_i}\}]^{2\theta-2} \left(-\bar{\alpha}\alpha e^{-2\lambda t_i} t_i \big/ \{1-\bar{\alpha}e^{-\lambda t_i}\}^2 - \alpha e^{-\lambda t_i} \big/ \{1-\bar{\alpha}e^{-\lambda t_i}\}\right)^2}{[1-[\alpha e^{-\lambda t_i} \big/ \{1-\bar{\alpha}e^{-\lambda t_i}\}]^{\theta}]^2}$$

$$+(1-m)\sum_{i=0}^{r}\frac{\theta(\theta-1)[\alpha e^{-\lambda t_i} \big/ \{1-\bar{\alpha}e^{-\lambda t_i}\}]^{\theta-2} \left(-\bar{\alpha}\alpha e^{-2\lambda t_i} t_i \big/ \{1-\bar{\alpha}e^{-\lambda t_i}\}^2 - \alpha e^{-\lambda t_i} \big/ \{1-\bar{\alpha}e^{-\lambda t_i}\}\right)^2}{1-[\alpha e^{-\lambda t_i} \big/ \{1-\bar{\alpha}e^{-\lambda t_i}\}]^{\theta}}$$

$$+(1-m)\sum_{i=0}^{r}\frac{\theta\, [\alpha e^{-\lambda t_i} \big/ \{1-\bar{\alpha}e^{-\lambda t_i}\}]^{\theta-1} \begin{pmatrix} -2\bar{\alpha}^2\alpha e^{-3\lambda t_i} t_i^2 \big/ \{1-\bar{\alpha}e^{-\lambda t_i}\}^3 + 3\bar{\alpha}\alpha e^{-2\lambda t_i} t_i^2 \big/ \{1-\bar{\alpha}e^{-\lambda t_i}\}^2 \\ +\alpha e^{-\lambda t_i} t_i^2 \big/ \{1-\bar{\alpha}e^{-\lambda t_i}\} \end{pmatrix}}{1-[\alpha e^{-\lambda t_i} \big/ \{1-\bar{\alpha}e^{-\lambda t_i}\}]^{\theta}}$$

$$\frac{\partial^2 \ell}{\partial m \partial n} = r\psi'(m+n)$$

$$\frac{\partial^2 \ell}{\partial m \partial \theta} = -\sum_{i=0}^{r}\frac{[\alpha e^{-\lambda t_i} \big/ \{1-\bar{\alpha}e^{-\lambda t_i}\}]^{\theta} \log[\alpha e^{-\lambda t_i} \big/ \{1-\bar{\alpha}e^{-\lambda t_i}\}]}{1-[\alpha e^{-\lambda t_i} \big/ \{1-\bar{\alpha}e^{-\lambda t_i}\}]^{\theta}}$$

$$\frac{\partial^2 \ell}{\partial m \partial \alpha} = -\sum_{i=0}^{r}\frac{\theta[\alpha e^{-\lambda t_i} \big/ \{1-\bar{\alpha}e^{-\lambda t_i}\}]^{\theta-1}\left(-[\alpha e^{-2\lambda t_i} \big/ \{1-\bar{\alpha}e^{-\lambda t_i}\}^2] + [\alpha e^{-\lambda t_i} \big/ \{1-\bar{\alpha}e^{-\lambda t_i}\}]\right)}{1-[\alpha e^{-\lambda t_i} \big/ \{1-\bar{\alpha}e^{-\lambda t_i}\}]^{\theta}}$$

$$\frac{\partial^2 \ell}{\partial m \partial \lambda} = -\sum_{i=0}^{r}\frac{\theta[\alpha e^{-\lambda t_i} \big/ \{1-\bar{\alpha}e^{-\lambda t_i}\}]^{\theta-1}\left(-[\bar{\alpha}\alpha e^{-2\lambda t_i} t_i \big/ \{1-\bar{\alpha}e^{-\lambda t_i}\}^2] - [\alpha e^{-\lambda t_i} t_i \big/ \{1-\bar{\alpha}e^{-\lambda t_i}\}]\right)}{1-[\alpha e^{-\lambda t_i} \big/ \{1-\bar{\alpha}e^{-\lambda t_i}\}]^{\theta}}$$

$$\frac{\partial^2 \ell}{\partial n \partial \theta} = \sum_{i=0}^{r}\log[\alpha e^{-\lambda t_i} t_i \big/ \{1-\bar{\alpha}e^{-\lambda t_i}\}]$$

$$\frac{\partial^2 \ell}{\partial n \partial \alpha} = \theta\sum_{i=0}^{r}\frac{e^{-\lambda t_i}\, [1-\bar{\alpha}e^{-\lambda t_i}]\left(-\alpha e^{-2\lambda t_i} \big/ \{1-\bar{\alpha}e^{-\lambda t_i}\}^2 + e^{-\lambda t_i} \big/ \{1-\bar{\alpha}e^{-\lambda t_i}\}\right)}{b}$$

$$\frac{\partial^2\ell}{\partial n\partial\lambda}=\theta\sum_{i=0}^{r}\frac{e^{-\lambda t_i}\,[1-\bar{\alpha}e^{-\lambda t_i}]\left(-\alpha\bar{\alpha}e^{-2\lambda t_i}t_i\Big/\{1-\bar{\alpha}e^{-\lambda t_i}\}^2+\alpha e^{-\lambda t_i}\,t_i\Big/\{1-\bar{\alpha}e^{-\lambda t_i}\}\right)}{b}$$

$$\frac{\partial^2\ell}{\partial\theta\partial\alpha}=\frac{r}{\alpha}+\sum_{i=0}^{r}\frac{e^{-\lambda t_i}}{1-\bar{\alpha}e^{-\lambda t_i}}$$

$$+(n-1)\sum_{i=0}^{r}\frac{e^{\lambda t_i}\,[1-\bar{\alpha}e^{-\lambda t_i}]\left(-\alpha\,e^{-2\lambda t_i}\Big/\{1-\bar{\alpha}e^{-\lambda t_i}\}^2+e^{-\lambda t_i}\Big/\{1-\bar{\alpha}e^{-\lambda t_i}\}\right)}{b}$$

$$+(1-m)\sum_{i=0}^{r}\frac{[\alpha e^{-\lambda t_i}\big/\{1-\bar{\alpha}e^{-\lambda t_i}\}]^{\theta\text{-}1}\left(-\alpha\,e^{-2\lambda t_i}\Big/\{1-\bar{\alpha}e^{-\lambda t_i}\}^2+e^{-\lambda t_i}\Big/\{1-\bar{\alpha}e^{-\lambda t_i}\}\right)}{1-[\alpha e^{-\lambda t_i}\big/\{1-\bar{\alpha}e^{-\lambda t_i}\}]^{\theta}}$$

$$-\theta(1\text{-}\mathrm{m})\sum_{i=0}^{r}\frac{[\alpha e^{-\lambda t_i}\big/\{1-\bar{\alpha}e^{-\lambda t_i}\}]^{2\theta\text{-}1}\left(-\alpha\,e^{-2\lambda t_i}\Big/\{1-\bar{\alpha}e^{-\lambda t_i}\}^2+e^{-\lambda t_i}\Big/\{1-\bar{\alpha}e^{-\lambda t_i}\}\right)\log[e^{-\lambda t_i}\big/\{1-\bar{\alpha}e^{-\lambda t_i}\}]}{[1-[\alpha e^{-\lambda t_i}\big/\{1-\bar{\alpha}e^{-\lambda t_i}\}]^{\theta}]^2}$$

$$-\theta(1\text{-}\mathrm{m})\sum_{i=0}^{r}\frac{[\alpha e^{-\lambda t_i}\big/\{1-\bar{\alpha}e^{-\lambda t_i}\}]^{\theta\text{-}1}\left(-\alpha\,e^{-2\lambda t_i}\Big/\{1-\bar{\alpha}e^{-\lambda t_i}\}^2+e^{-\lambda t_i}\Big/\{1-\bar{\alpha}e^{-\lambda t_i}\}\right)\log[e^{-\lambda t_i}\big/\{1-\bar{\alpha}e^{-\lambda t_i}\}]}{1-[\alpha e^{-\lambda t_i}\big/\{1-\bar{\alpha}e^{-\lambda t_i}\}]^{\theta}}$$

$$\frac{\partial^2\ell}{\partial\theta\partial\lambda}=-\sum_{i=0}^{r}t_i+\sum_{i=0}^{r}\log[\bar{\alpha}e^{-\lambda t_i}\,t_i\big/\{1-\bar{\alpha}e^{-\lambda t_i}\}]$$

$$+(n-1)\sum_{i=0}^{r}\frac{e^{-\lambda t_i}\,[1-\bar{\alpha}e^{-\lambda t_i}]\left(-\alpha\bar{\alpha}e^{-2\lambda t_i}t_i\Big/\{1-\bar{\alpha}e^{-\lambda t_i}\}^2-\alpha e^{-\lambda t_i}\,t_i\Big/\{1-\bar{\alpha}e^{-\lambda t_i}\}\right)}{b}$$

$$+(1-m)\sum_{i=0}^{r}\frac{[\alpha e^{-\lambda t_i}\big/\{1-\bar{\alpha}e^{-\lambda t_i}\}]^{\theta\text{-}1}\left(-\bar{\alpha}\alpha\,e^{-2\lambda t_i}t_i\Big/\{1-\bar{\alpha}e^{-\lambda t_i}\}^2-e^{-\lambda t_i}\Big/\{1-\bar{\alpha}e^{-\lambda t_i}\}\right)}{1-[\alpha e^{-\lambda t_i}\big/\{1-\bar{\alpha}e^{-\lambda t_i}\}]^{\theta}}$$

$$-\theta(1\text{-}\mathrm{m})\sum_{i=0}^{r}\frac{[\alpha e^{-\lambda t_i}\big/\{1-\bar{\alpha}e^{-\lambda t_i}\}]^{2\theta\text{-}1}\left(-\alpha\bar{\alpha}e^{-2\lambda t_i}\Big/\{1-\bar{\alpha}e^{-\lambda t_i}t_i\}^2+e^{-\lambda t_i}t_i\Big/\{1-\bar{\alpha}e^{-\lambda t_i}\}\right)\log[e^{-\lambda t_i}\big/\{1-\bar{\alpha}e^{-\lambda t_i}\}]}{[1-[\alpha e^{-\lambda t_i}\big/\{1-\bar{\alpha}e^{-\lambda t_i}\}]^{\theta}]^2}$$

$$-\theta(1\text{-}\mathrm{m})\sum_{i=0}^{r}\frac{[\alpha e^{-\lambda t_i}\big/\{1-\bar{\alpha}e^{-\lambda t_i}\}]^{\theta\text{-}1}\left(-\alpha\bar{\alpha}e^{-2\lambda t_i}t_i\Big/\{1-\bar{\alpha}e^{-\lambda t_i}\}^2+e^{-\lambda t_i}t_i\Big/\{1-\bar{\alpha}e^{-\lambda t_i}\}\right)\log[e^{-\lambda t_i}\big/\{1-\bar{\alpha}e^{-\lambda t_i}\}]}{1-[\alpha e^{-\lambda t_i}\big/\{1-\bar{\alpha}e^{-\lambda t_i}\}]^{\theta}}$$

$$\frac{\partial^2\ell}{\partial\alpha\partial\lambda}=(\theta+1)\sum_{i=0}^{r}\left(-\frac{\bar{\alpha}e^{-2\lambda t_i}\,t_i}{[1-\bar{\alpha}e^{-\lambda t_i}]^2}-\frac{e^{-\lambda t_i}t_i}{1-\bar{\alpha}e^{-\lambda t_i}}\right)$$

$$+\theta(n-1)\sum_{i=0}^{r}\frac{t_i\left(-\alpha e^{-2\lambda t_i}\Big/\{1-\bar{\alpha}e^{-\lambda t_i}\}^2+e^{-\lambda t_i}\Big/\{1-\bar{\alpha}e^{-\lambda t_i}\}\right)}{b}$$

$$+\theta(n-1)\sum_{i=0}^{r}\frac{e^{-\lambda t_i}\,[1-\bar{\alpha}e^{-\lambda t_i}]t_i\left(-\alpha e^{-2\lambda t_i}\Big/\{1-\bar{\alpha}e^{-\lambda t_i}\}^2+e^{-\lambda t_i}\Big/\{1-\bar{\alpha}e^{-\lambda t_i}\}\right)}{b}$$

$$+\theta(n-1)\sum_{i=0}^{r}\frac{e^{-\lambda t_i}[1-\bar{\alpha}e^{-\lambda t_i}]\begin{pmatrix}2\bar{\alpha}\alpha e^{-3\lambda t_i}t_i/\{1-\bar{\alpha}e^{-\lambda t_i}\}^3-\bar{\alpha}e^{-2\lambda t_i}t_i/\{1-\bar{\alpha}e^{-\lambda t_i}\}^2\\+2\alpha e^{-2\lambda t_i}t_i/\{1-\bar{\alpha}e^{-\lambda t_i}\}^2-e^{-\lambda t_i}t_i/\{1-\bar{\alpha}e^{-\lambda t_i}\}\end{pmatrix}}{b}$$

$$+(1-m)\sum_{i=0}^{r}\frac{\theta\{\alpha e^{-\lambda t_i}/(1-\bar{\alpha}e^{-\lambda t_i})\}^{\theta-1}\begin{pmatrix}2\bar{\alpha}\alpha e^{-3\lambda t_i}t_i/\{1-\bar{\alpha}e^{-\lambda t_i}\}^3-\bar{\alpha}e^{-2\lambda t_i}t_i/\{1-\bar{\alpha}e^{-\lambda t_i}\}^2\\+2\alpha e^{-2\lambda t_i}t_i/\{1-\bar{\alpha}e^{-\lambda t_i}\}^2-e^{-\lambda t_i}t_i/\{1-\bar{\alpha}e^{-\lambda t_i}\}\end{pmatrix}}{1-\{\alpha e^{-\lambda t_i}/(1-\bar{\alpha}e^{-\lambda t_i})\}^{\theta}}$$

$$+(1-m)\sum_{i=0}^{r}\frac{\begin{matrix}\theta^2\{\alpha e^{-\lambda t_i}/(1-\bar{\alpha}e^{-\lambda t_i})\}^{2\theta-2}(-\alpha e^{-2\lambda t_i}/\{1-\bar{\alpha}e^{-\lambda t_i}\}^2+e^{-\lambda t_i}/\{1-\bar{\alpha}e^{-\lambda t_i}\})\\ \times\left(-(\bar{\alpha}\alpha e^{-2\lambda t_i}t_i/\{1-\bar{\alpha}e^{-\lambda t_i}\}^2-\alpha e^{-\lambda t_i}t_i/\{1-\bar{\alpha}e^{-\lambda t_i}\})\right)\end{matrix}}{[1-\{\alpha e^{-\lambda t_i}/(1-\bar{\alpha}e^{-\lambda t_i})\}^{\theta}]^2}$$

$$+(1-m)\sum_{i=0}^{r}\frac{\begin{matrix}\theta(\theta-1)\{\alpha e^{-\lambda t_i}/(1-\bar{\alpha}e^{-\lambda t_i})\}^{\theta-2}(-\alpha e^{-2\lambda t_i}/\{1-\bar{\alpha}e^{-\lambda t_i}\}^2+e^{-\lambda t_i}/\{1-\bar{\alpha}e^{-\lambda t_i}\})\\ \times\left(-(\bar{\alpha}\alpha e^{-2\lambda t_i}t_i/\{1-\bar{\alpha}e^{-\lambda t_i}\}^2-\alpha e^{-\lambda t_i}t_i/\{1-\bar{\alpha}e^{-\lambda t_i}\})\right)\end{matrix}}{1-\{\alpha e^{-\lambda t_i}/(1-\bar{\alpha}e^{-\lambda t_i})\}^{\theta}}$$

Where $\psi'(.)$ is the derivative of the digamma function.